\documentclass[11pt]{amsart}

\usepackage{geometry}
\usepackage{amsmath, amsthm, amsfonts, hyperref, graphicx, ifpdf, mathrsfs,color}
\usepackage{amssymb}
\usepackage{faktor}
\usepackage{bm}
\usepackage{enumerate}
\usepackage[utf8]{inputenc}
\usepackage{mathtools}
\usepackage{esint}
\usepackage{tikz-cd} 
\usepackage{verbatim}
\usepackage{setspace}

\hypersetup{
	unicode=false,          
	pdftoolbar=true,        
	pdfmenubar=true,        
	pdffitwindow=false,     
	pdfstartview={FitH},    
	pdftitle={My title},    
	pdfauthor={Author},     
	pdfsubject={Subject},   
	pdfcreator={Creator},   
	pdfproducer={Producer}, 
	pdfkeywords={keywords}, 
	pdfnewwindow=true,      
	colorlinks=true,       
	linkcolor=blue,          
	citecolor=blue,        
	filecolor=magenta,      
	urlcolor=MidnightBlue          
}

\theoremstyle{plain}
\newtheorem{theorem}{Theorem}[section]
\newtheorem{lem}[theorem]{Lemma}
\newtheorem{pro}[theorem]{Proposition}
\newtheorem{cor}[theorem]{Corollary}

\theoremstyle{definition}
\newtheorem{definition}[theorem]{Definition}

\theoremstyle{remark}
\newtheorem{rem}[theorem]{Remark}

\numberwithin{equation}{section}

\newcommand{\afm}{almost-Fuchsian manifold}

\newcommand{\app}{asymptotic Plateau problem}

\newcommand{\cvh}{convex hull}

\newcommand{\Hd}{Hausdorff dimension}

\newcommand{\htm}{hyperbolic three-manifold}

\newcommand{\Jc}{Jordan curve}

\newcommand{\lmt}{limit set}

\newcommand{\maxp}{maximum principle}

\newcommand{\MS}{moduli space}
\newcommand{\ms}{minimal surface}
\newcommand{\md}{minimal disk}

\newcommand{\pc}{principal curvature}
\newcommand{\qf}{quasi-Fuchsian}
\newcommand{\qfg}{quasi-Fuchsian group}
\newcommand{\qfm}{quasi-Fuchsian manifold}

\newcommand{\spc}{small curvature}

\newcommand{\TS}{Teichm\"{u}ller space}

\newcommand{\tm}{three-manifold}
\newcommand{\tg}{totally geodesic}

\newcommand{\be}{\begin{equation}}
	\newcommand{\ene}{\end{equation}}
\newcommand{\br}{\begin{rem}}
	\newcommand{\er}{\end{rem}}
\newcommand{\bl}{\begin{lem}}
	\newcommand{\bcor}{\begin{cor}}
		\newcommand{\ecor}{\end{cor}}
	\newcommand{\el}{\end{lem}}
\newcommand{\bd}{\begin{Def}}
	\newcommand{\ed}{\end{Def}}
\newcommand{\ben}{\begin{enumerate}}
	\newcommand{\een}{\end{enumerate}}
\newcommand{\bp}{\begin{proof}}
	\newcommand{\ep}{\end{proof}}
\newcommand{\bpo}{\begin{pro}}
	\newcommand{\epo}{\end{pro}}
\newcommand{\beq}{\begin{equation*}}
	\newcommand{\eeq}{\end{equation*}}
\newcommand{\bear}{\begin{eqnarray*}}
	\newcommand{\eear}{\end{eqnarray*}}
\newcommand{\bt}{\begin{theorem}}
	\newcommand{\et}{\end{theorem}}

\newcommand{\Sph}{\mathbb{S}}
\newcommand{\HH}{\mathbb{H}^3}

\newcommand{\HHnp}{\mathbb{H}^{n+1}}

\newcommand{\R}{\mathbb{R}}
\newcommand{\C}{\mathbb{C}}

\newcommand{\N}{\mathbb{N}}

\newcommand{\Cal}{\mathcal{C}}

\newcommand{\Z}{\mathbb{Z}}

\newcommand{\Cat}{\mathrm{Cat}}
\newcommand{\SCat}{\mathrm{Cat}^{\mathrm{solid}}}

\numberwithin{equation}{section}

\allowdisplaybreaks

\def\XXint#1#2#3{{\setbox0=\hbox{$#1{#2#3}{\int}$}
		\vcenter{\hbox{$#2#3$}}\kern-.5\wd0}}

\makeatletter
\def\@citestyle{\m@th\upshape\mdseries}
\def\citeform#1{{\bfseries#1}}
\def\@cite#1#2{{%
		\@citestyle[\citeform{#1}\if@tempswa, #2\fi]}}
\@ifundefined{cite }{%
	\expandafter\let\csname cite \endcsname\cite
	\edef\cite{\@nx\protect\@xp\@nx\csname cite \endcsname}%
}{}
\makeatother

\begin{document}
	\parskip1ex
	\title[Uniqueness and non-uniqueness]
	{Uniqueness and non-uniqueness for the asymptotic Plateau problem in hyperbolic space}

	\author{Zheng Huang}
	\address[Z. ~H.]{Department of Mathematics, The City University of New York, Staten Island, NY 10314, USA}
	\address{The Graduate Center, The City University of New York, 365 Fifth Ave., New York, NY 10016, USA}
	\email{zheng.huang@csi.cuny.edu}
	
	\author{Ben Lowe}
	\address[B.~L.]{Department of Mathematics, University of Chicago, Chicago, IL 60637, USA }
	\email{loweb24@uchicago.edu}
	
	\author{Andrea Seppi}
	\address[A.~S.]{Univ. Grenoble Alpes, CNRS, IF, 38000 Grenoble, France}
	\email{andrea.seppi@univ-grenoble-alpes.fr}

	\date{\today}
	
	\subjclass[2020]{Primary 53C42, 57K32}
	
	\begin{abstract}
		We prove several results on the number of solutions to the asymptotic Plateau problem in $\mathbb H^3$. Firstly we discuss criteria that ensure uniqueness. Given a Jordan curve $\Lambda$ in the asymptotic boundary of $\mathbb H^3$, we show that uniqueness of the minimal surfaces with asymptotic boundary $\Lambda$ is equivalent to uniqueness in the smaller class of stable minimal disks. 
		
		Then we show that if a quasicircle (or more generally, a Jordan curve of finite width) $\Lambda$ is the asymptotic boundary of a minimal surface $\Sigma$ with principal curvatures less than or equal to 1 in absolute value, then uniqueness holds.

		In the direction of non-uniqueness, we construct an example of a quasicircle that is the asymptotic boundary of uncountably many pairwise distinct stable minimal disks. 
		
	\end{abstract}
	
	\maketitle

	\vspace{-1cm}
	
	\begin{singlespace}
		\tableofcontents  
	\end{singlespace}

	\section {Introduction}

	The classical ``asymptotic Plateau problem" asks, given a {\Jc} $\Lambda$ on ${\Sph}^2_{\infty} = \partial_\infty \HH$, the existence and multiplicity of (properly embedded) {\ms}s $\Sigma$ in $\HH$, if any, that are asymptotic to $\Lambda$, in the sense that the closure of $\Sigma$ in ${\Sph}^2_{\infty}\cup \HH$ is equal to $\Lambda \cup \Sigma$. The existence of {\md} solutions to the asymptotic Plateau problem was first stated by Anderson (\cite{And83}). Using geometric measure theory, Anderson also obtained existence results in higher dimensions. 
	
	The uniqueness does not hold in general: as shown in \cite{And83, HW15}, taking advantage of group actions, one can construct a {\Jc} $\Lambda$ in ${\Sph}^2_{\infty} = \partial_\infty \HH$ which is the {\lmt} of 
	some {\qf} group such that $\Lambda$ spans multiple  {\md}s (even an arbitrarily large, but finite, number). Anderson (\cite{And83}) even constructed a curve $\Lambda$ which 
	spans infinitely many {\ms}s (the surfaces he constructs have positive genus).  On the other hand, when $\Lambda$ 
	is a round circle, the unique minimal surface it spans is a {\tg} disk. To look for unique solutions, it is therefore natural to consider the class of {\ms}s that are ``close" to {\tg}, for which $\Lambda$ is ``close" to a round circle. Related questions with conditions on natural invariants of $\Lambda$ were studied in \cite{Sep16} 
	(for the quasi-conformal constant of $\Lambda$), and \cite{HW13, San18} (for the {\Hd} of $\Lambda$). Some properness questions for the asymptotic Plateau problem solutions for various classes of curves were addressed, for example, in  \cite{GS00, AM10, co14}. 
	
	Motivated by the above results and by many natural questions arising from the study of the {\app}, in this paper we address two basic questions, which, for the sake of simplicity, we state only in dimension three, but may naturally be extended to hypersurfaces in $\mathbb{H}^{n+1}$:
	\begin{itemize}
		\item Under what conditions does a {\Jc} $\Lambda$ on ${\Sph}^2_{\infty}$ span exactly one {\ms} in $\HH$? 
		\item Does 
		there exist a {\Jc} $\Lambda$ on ${\Sph}^2_{\infty}$ that spans   infinitely many {\md}s in $\HH$, and if yes, which cardinality may the set of solutions have?   
	\end{itemize}

	\subsection{Characterizing uniqueness}
	
	Recall that a (hyper)surface is \emph{minimal}  if it is a critical point of the area under compactly supported variations. It is \emph{stable} if moreover the second variation of the area under any compactly supported variation is non-negative. In this introduction, we always implicitly consider \emph{properly embedded} (hyper)surfaces. Our first theorem shows that it suffices to check uniqueness in the class of \emph{stable minimal disks}. 
	
	\bt\label{thm:char uniqueness}
	Let $\Lambda$ be a {\Jc} on ${\Sph}^2_{\infty} = \partial_{\infty} \HH$. Then $\Lambda$ spans a unique  minimal surface if and only if it spans a unique stable minimal disk.
	\et
	
	A statement similar to Theorem \ref{thm:char uniqueness}, but in the context of the \emph{finite} Plateau problem, was proved in \cite{MY19}. When $\Lambda$ is the limit set of a quasi-Fuchsian group, applying the results of \cite{GLP21}, we can prove a stronger statement.
	
	\bt\label{thm:char uniqueness invariant}
	Let $\Lambda\subset{\Sph}^2_{\infty} = \partial_{\infty} \HH$ be the limit set of a quasi-Fuchsian group $\Gamma<\mathrm{Isom}(\HH)$ isomorphic to the fundamental group of a closed surface. Then $\Lambda$ spans a unique  minimal surface if and only if it spans a unique $\Gamma$-invariant stable minimal disk.
	\et
	
	
	The main idea in the proof of Theorem \ref{thm:char uniqueness} is the following: we show that if there is a {\ms} $\Sigma$, which is not stable or is not topologically a disk, with asymptotic boundary 
	$\Lambda$, then we can construct two distinct --- actually, disjoint --- stable minimal disks with the same asymptotic boundary $\Lambda$ (Theorem \ref{main3}). These two stable minimal disks are constructed ``on each side'' of $\Sigma$.
	
	Based on Theorem \ref{thm:char uniqueness}, to prove Theorem \ref{thm:char uniqueness invariant} it then suffices to prove that, if $\Lambda$ spans a unique $\Gamma$-invariant stable minimal disk, then it spans a unique stable minimal disk. Then we apply one of the results in \cite{GLP21}, which provides, for any quasi-Fuchsian manifold, the existence of a foliation where the sign of the mean curvature is constant. In particular, if the $\Gamma$-invariant stable minimal disk is unique, then the leaves of the foliation have positive mean curvature on one side of the minimal surface, and negative mean curvature on the other side.  Lifting this foliation to a foliation of $\HH$ is then the essential ingredient for our proof of Theorem \ref{thm:char uniqueness invariant}, together with well-known uniform bounds on the second fundamental form of stable minimal surfaces.
	
	
	
	\subsection{Uniqueness criteria via curvature conditions}
	
	Next, we turn our attention to sufficient conditions for uniqueness. For an immersed hypersurface in $\HHnp$, or more generally in a hyperbolic ($n+1$)-manifold, we say it has {\it strongly small curvature} if its {\pc}s $\{\lambda_i\}$ satisfy that 
	\be\label{cc}
	|\lambda_i| \le 1-\epsilon, \ i=1,\ldots,n, \ \ \ \text{for some small} \ \ \epsilon>0. 
	\ene
	Similarly it has {\it small curvature} if $|\lambda_i| < 1$, and it has {\it weakly small curvature} if $|\lambda_i| \le 1$. This definition has some immediate consequences: for instance, a complete immersion of weakly small curvature is in fact a properly embedded disk (\cite{Eps84, Eps86}, see also \cite{ES22} for a generalization). See Subsection \ref{subsec:small curvatures} and Appendix \ref{appendix} for more details.

	Surfaces of small curvature are very special in {\tm} theory: Thurston observed that a closed surface of {\spc} in a complete {\htm}	is incompressible (\cite{Thu86, Lei06}); they are abundant in closed {\htm}s (\cite{KM12}); many results have been extended to the study of complete noncompact {\htm} of finite volume (\cite{Rub05, CF19, kw21}). It is often favorable to consider canonical representatives within a homotopy class of surfaces, and minimal surfaces are in many ways the most natural choice. Among hyperbolic three-manifolds, \emph{almost-Fuchsian} manifolds are quasi-Fuchsian manifolds which admit a closed {\ms} of {\spc}. 
	This notion was introduced by Uhlenbeck and it played an important role in her study of parametrization of the {\MS} of {\ms}s in {\htm}s (\cite{Uhl83}). Subsequently many different aspects of this subclass of {\qfm} have been studied up to recent years (for instance \cite{KS07, GHW10, HL21} and many others). 
	
	It is known (\cite{Uhl83}) that any {\afm} admits a unique closed {\ms} --- in other words, identifying the almost-Fuchsian manifold with a quotient $\HH/\Gamma$, the limit set $\Lambda$ of the group $\Gamma$ bounds a unique $\Gamma$-invariant   minimal disk asymptotic to $\Lambda$.
	Inspired by this fact, we prove (Corollary \ref{cor strongly small} below) that if a {\Jc} $\Lambda$ (not necessarily group equivariant) spans a {\md} $\Sigma$ of {\underline{strongly} {\spc}} in $\HH$, then $\Sigma$ is the unique {\ms} asymptotic to $\Lambda$. Our results, however, are more general. The main result we prove in this direction is the following:


	\bt\label{main1}
	Let $\Lambda$ be a Jordan curve on ${\Sph}^2_{\infty} = \partial_{\infty} \HH$ of finite width, and let $\Sigma$ be a  minimal surface in $\HH$ of 
	{\underline {weakly}} {\spc} asymptotic to $\Lambda$. Then $\Sigma$ is the unique {\ms} in $\HH$ asymptotic to 
	$\Lambda$. Moreover, $\Sigma$ is area-minimizing.
	\et 
	
	
	Let us explain the terminology of the statement. First, the \emph{width} of a Jordan curve $\Lambda$ in ${\Sph}^2_{\infty} = \partial_{\infty} \HH$ has been introduced in \cite{BDMS21} as the supremum over all points in the convex hull of $\Lambda$ of the sum of the distances from each boundary component of the convex hull. Quasicircles are an important class of Jordan curves, which are known to have finite width. (On the other hand, there exist Jordan curves of finite width which are not quasicircles, as constructed in \cite{BDMS21}.) Hence the conclusion of Theorem \ref{main1} in particular holds for any quasicircle $\Lambda$.
	
	Second, we recall that a hypersurface $\Sigma$ is \emph{area-minimizing} if any compact codimension-zero submanifold with boundary has smaller area than any rectifiable hypersurface with the same boundary in the ambient space. This implies that $\Sigma$ is a stable minimal hypersurface. 
	
	We remark that the setting of Theorems \ref{thm:char uniqueness} and \ref{thm:char uniqueness invariant} is more general than Theorem \ref{main1}. Indeed, it follows from \cite[Theorem 5.2]{HL21} that there are examples of quasi-Fuchsian groups $\Gamma$ whose limit set $\Lambda$ bounds a unique $\Gamma$-invariant stable minimal disk $\Sigma$ (hence, by Theorem \ref{thm:char uniqueness invariant}, a unique minimal surface), but $\Sigma$ does \textit{not} have weakly small curvature.
	
	Let us sketch the proof of Theorem \ref{main1}, which is based on elementary geometric arguments. The first observation is that the finite width condition implies that every minimal surface $\Sigma'$ has, roughly speaking, finite normal distance from $\Sigma$. Moreover, by Theorem \ref{thm:char uniqueness}, it suffices to show uniqueness among \emph{stable} minimal disks. Now, if the maximum normal distance between $\Sigma$ and another stable {\ms} is realized in the interior, the conclusion follows from a {\maxp} argument, taking advantage  the properties of the normal flow from $\Sigma$. If the maximum normal distance is not realized in the interior, we use isometries to send the run-away sequence of points back to a fixed point and use the compactness theorems for stable {\ms}s, so as to reduce essentially to the previous case.

	We now derive a consequence of Theorem \ref{main1}. Observe that if a properly embedded (hyper)surface $\Sigma$ has {\underline{strongly}} small curvatures, then its asymptotic boundary has finite width (see Lemma \ref{lemma finite width} in Appendix \ref{appendix}). Hence we obtain: 
	
	\bcor\label{cor strongly small}
	Let $\Lambda$ be a Jordan curve on ${\Sph}^2_{\infty} = \partial_{\infty} \HH$, and let $\Sigma$ be a minimal surface in $\HH$ of 
	{\underline {strongly}} {\spc} asymptotic to $\Lambda$. Then $\Sigma$ is the unique   minimal surface in $\HH$ asymptotic to 
	$\Lambda$. Moreover, $\Sigma$ is area-minimizing.
	\ecor

	In Section \ref{sec:higher} we do a brief excursion in higher dimensions, and show that Corollary \ref{cor strongly small} actually holds  in $\HHnp$ (Theorem \ref{thm strongly small higher dim}). That is, if $\Lambda$ is a topologically embedded $(n-1)$-sphere on ${\Sph}^n_{\infty} = \partial_{\infty} \HHnp$ that spans a minimal hypersurface in $\HHnp$ of 
	{{strongly}} {\spc}, then $\Sigma$ is the unique   minimal hypersurface in $\HHnp$ asymptotic to 
	$\Lambda$ (and is area-minimizing).
	The proof rests on an application of a general maximum principle proved in \cite{Whi10} in the context of minimal varifolds (see also \cite{jt03}). This approach, however, does not lead to the analogue of Theorem \ref{main1} for minimal hypersurfaces in any dimensions, which would be a stronger statement.
	
	


	\subsection{Strong non-uniqueness}
	
	Now we turn to the other extreme case: we aim to construct a {\Jc} on $\Sph^2_\infty$ which spans \emph{a lot} of {{\md}s}. To give some context, Hass-Thurston conjectured that no closed hyperbolic 3-manifold admits a foliation by minimal surfaces.  Anderson even conjectured (\cite{And83}) that no hyperbolic 3-manifold admits a local 1-parameter family of closed minimal surfaces, and proved this statement for quasi-Fuchsian hyperbolic 3-manifolds. It has been a folklore conjecture that no Jordan curve in ${\Sph}^2_{\infty}= \partial_\infty \mathbb{H}^3$ asymptotically bounds a 1-parameter family of minimal surfaces. The full extent of these conjectures remain as major open questions in the field.
	
	What has been proven up to this point tends to support these conjectures.  Huang-Wang \cite{HW19} and Hass \cite{has15} made progress on the Hass-Thurston conjecture for certain fibered closed hyperbolic 3-manifolds containing short geodesics; Wolf-Wu (\cite{ww20}) ruled out so-called geometric local 1-parameter families of closed minimal surfaces; it follows from the work of Alexakis-Mazzeo \cite{AM10} that a generic $C^{3,\alpha}$ simple closed curve in the boundary at infinity of $\mathbb{H}^3$ bounds only finitely many minimal surfaces of any given finite genus; Coskunuzer proved that generic simple closed curves in $\partial_\infty \mathbb{H}^3$ bound unique area-minimizing surfaces \cite{c11genericuniqueness}.  
	
	What we prove, while compatible with the folklore conjecture, is in the other direction. Based on the aforementioned results, one might be tempted to strengthen the folklore conjecture to the statement that any Jordan curve in ${\Sph}^2_{\infty}= \partial_\infty \mathbb{H}^3$ bounds at most countably many minimal surfaces.  We show that this stronger statement is false:

	\bt\label{thm:uncountably}
	There exists a quasicircle in ${\Sph}^2_{\infty} = \partial_{\infty} \HH$ spanning uncountably many pairwisely distinct stable minimal disks.
	\et
	Let us emphasize some important features of the construction of this extreme curve $\Lambda$. In \cite{And83}, Anderson constructed a {\Jc} 
	which is the limit set of a {\qfg} (hence a ``fractal'' Jordan curve with Hausdorff dimension strictly greater than one) such that it spans infinitely many {\ms}s, one of which is a 
	{\md}. In \cite{HW15}, for each integer $N > 1$, also using the limit set of a {\qfg}, an extreme curve spanning at least 
	$2^N$ distinct {\md}s was constructed. However, Anderson (\cite{And83}) has shown that any {\qfm} only admits finitely many least area closed {\ms}s 
	diffeomorphic to the fiber, which poses a possible limitation on how much one can improve the aforementioned constructions to find infinitely many 
	{\md}s if one insists on using the limit set of some {\qfg} as the curve at infinity. The starting point of our  construction is similar to the ideas in \cite{HW15}, but the Jordan curve is constructed in such a way to allow an improvement of the argument, leading to $2^{\mathbb N}$ pairwise distinct minimal disk. 
	Moreover, since this Jordan curve is not invariant under any {\qfg}, we must adopt a different approach in order to produce the minimal disks, in a spirit similar to the proof of Theorem \ref{thm:char uniqueness}.


	\subsection{Quasiconformal constant}
	
	We conclude this introduction with an improvement of the curvature estimates obtained in \cite{Sep16} in terms of quasiconformal constants, by a direct application of Theorem \ref{main1}. More concretely, \cite[Theorem A]{Sep16} showed that there exist universal constants $C>0$ and $K_0>1$ such that any stable minimal disk in $\HH$ with asymptotic boundary a $K$-quasicircle, for $K<K_0$, has principal curvatures bounded in absolute value by $C\log K$. This result has been recently applied in several directions, see \cite{Bis19,Low21,CMN22,KMS23}.
	
	The proof, however, relies on the application of compactness for minimal surfaces, and therefore requires stability. However, when $K$ is sufficiently small, the principal curvatures of the area-minimizing (hence stable) disk whose existence is guaranteed by \cite{And83} are less than $1-\epsilon$ in absolute value, and therefore, as a consequence of Theorem \ref{main1}, $\Sigma$ is the unique minimal surface. Up to taking a smaller constant $K_0$, we can therefore remove the stability assumption:
	
	\bcor\label{qc}
	There exist universal constants $C>0$ and $K_0>1$ such that the principal curvatures $\lambda_i$ of any minimal surface $\Sigma$ in $\HH$ with asymptotic boundary a $K$-quasicircle with $K\leq K_0$ satisfy 
	$$|\lambda_i|\leq C\log K,\; i=1,2.$$
	\ecor
	
	In particular, we also improve \cite[Theorem B]{Sep16} (up to choosing a smaller constant) by removing the stability assumption. 
	
	\bcor\label{qc2}
	There exists a universal constant $K'_0>1$ such that any $K$-quasicircle with $K\leq K_0'$ is the asymptotic boundary of a unique minimal surface, which is an area-minimizing disk of strongly small curvature.
	\ecor
	\subsection{Organization of the paper} In the preliminary section \S \ref{prelim}, we collect and prove some facts that will be useful in proving our main results introduced above. In \S \ref{sec:APP}, we give a complete proof of the existence theorem for a disk solution to the {\app}, initially stated in \cite{And83}. In fact a variation of our approach to the {\app}, in a slightly more complicated setting, will be an essential ingrendient for the proofs of Theorems \ref{thm:char uniqueness} and \ref{thm:uncountably}. In \S \ref{sec3} we prove Theorems \ref{thm:char uniqueness}, \ref{thm:char uniqueness invariant} and Theorem \ref{main1}. In \S \ref{sec:higher} we prove a generalization of Corollary \ref{cor strongly small} in higher dimensions, and we briefly outline some known results in higher codimension. In \S \ref{sec:uncountably} we detail a construction of an extreme {\Jc} which spans uncountably many {\md}s in $\HH$ and hence prove Theorem \ref{thm:uncountably}. In Appendix \ref{appendix} we provide the details to extend some well-known arguments for hypersurfaces of {\spc} to the setting of weakly {\spc} that is of interest here. 
	\subsection{Acknowledgements}
	We would like to thank Baris Coskunuzer for pointing out the gaps in the paper (\cite{And83}), Bill Meeks for his insightful suggestions and for answering questions on minimal laminations, Franco Vargas Pallete for comments on an earlier version of this manuscript, and Biao Wang for his generous help. 
	
	
	Part of this work was done during a visit of the first-named author at the Institut Fourier (Université Grenoble Alpes) in the framework of the ``Visiting Scientist Campaign 2023'', he wishes to thank the institute for excellent working environment. The second-named author was supported by an NSF grant DMS-2202830. The third-named author was funded by the ANR JCJC grant GAPR (ANR-22-CE40-0001) and by the European Union (ERC, GENERATE, 101124349). Views and opinions expressed are however those of the author(s) only and do not necessarily reflect those of the European Union or the European Research Council Executive Agency. Neither the European Union nor the granting authority can be held responsible for them.
	
	
	\section{Preliminaries} \label{prelim} 
	
	In this article, we denote by $\HHnp$ the hyperbolic space of dimension $n+1$, and by $\Sph^n_\infty=\partial_\infty\HHnp$ its visual boundary. 
	
	\subsection{Hypersurfaces theory} 
	
	Given an immersed hypersurface $\Sigma$ in $\HHnp$, we recall that its first fundamental form is the restriction to $T\Sigma$ of the  hyperbolic metric $h$ of $\HHnp$. The second fundamental form is defined as $$A_\Sigma(v,w)=h(\nabla^{h}_{v}W,N_\Sigma)~,$$
	where $\nabla^{h}$ is the Levi-Civita connection of $h$, $W$ is a local smooth extension of $w$, and  $N_\Sigma:\Sigma\to T\HHnp$ is a continuous choice of a unit normal vector to the immersion. (Since $N_\Sigma$ is uniquely determined only up to a sign, so too is $A_\Sigma$.) 
	
	The second fundamental form satisfies the identity
	$$A_\Sigma(v,w)=h(B_\Sigma(v),w)$$
	where $B_\Sigma$ is the shape operator, namely the endomorphism of the tangent bundle of $\Sigma$ given by $B_\Sigma(v)=-\nabla^h_vN_\Sigma$. The principal curvatures of $\Sigma$ are the eigenvalues of $B_\Sigma$, denoted by $\lambda_1,\ldots,\lambda_n$. 
	
	The mean curvature of the immersed hypersurface $\Sigma$ is 
	$$H_\Sigma=\mathrm{tr}(B_\Sigma)=\lambda_1+\ldots+\lambda_n~,$$
	and $\Sigma$ is  minimal if and only if its mean curvature vanishes identically. Although $H_\Sigma$ depends (but only up to a sign) on the choice of the normal vector $N_\Sigma$, the condition of being minimal does not. Also, the mean curvature vector, which is defined as $H_\Sigma\cdot N_\Sigma$, does not depend on such a choice.
	
	Finally, the norm of the second fundamental form is 
	$$\|A_\Sigma\|=\sqrt{\mathrm{tr}(B_\Sigma B_\Sigma^T)}=\sqrt{\lambda_1^2+\ldots+\lambda_n^2}~.$$

	\subsection{Convex hull and width of a {\Jc}} 
	Let  $\Lambda$ be a topologically embedded $n$-sphere in $\Sph^n_\infty =\partial\HHnp$. We recall that the {\cvh} 
	$\Cal(\Lambda) \subset \overline{\HHnp}=\HHnp\cup\Sph^n_\infty$ of $\Lambda$ is the smallest geodesically convex subset that contains $\Lambda$. 
	When $\Lambda$ is not the boundary of a totally geodesic hyperplane in $\HHnp$, $\Cal(\Lambda)$ is homeomorphic to a ball, and, by the Jordan-Brouwer separation theorem, its boundary 
	is the union of $\Lambda$ and two properly embedded disks, denoted by $\partial^{+}\Cal(\Lambda)$ and $\partial^{-}\Cal(\Lambda)$. 
	When $\Lambda$ is the boundary of a totally geodesic hyperplane $P$, $\Cal(\Lambda)$ equals $P\cup\Lambda$. In this case, by an abuse of notation, we will still use the symbols $\partial^{+}\Cal(\Lambda)$ and $\partial^{-}\Cal(\Lambda)$, meaning that $P=\partial^{+}\Cal(\Lambda)=\partial^{-}\Cal(\Lambda)$.
	

	Following (and extending to all dimensions) the recent work \cite{BDMS21}, we now define the {\it width} of $\Lambda$. 
	
	\begin{definition}\label{width}
		Given a topologically embedded $n$-sphere $\Lambda$ in $\Sph_\infty^n$, the width of $\Lambda$ is defined as:
		$$w(\Lambda)=\sup_{x\in\Cal(\Lambda)}\left(d(x,\partial^{+}\Cal(\Lambda))+d(x,\partial^{-}\Cal(\Lambda))\right)\in[0,+\infty]~.$$
	\end{definition}

	The following lemma, that relates minimal hypersurfaces and the convex hull of their asymptotic boundaries,  is well-known.
	
	\bl\label{lemma interior convex hull}
	Given a topologically embedded $(n-1)$-sphere $\Lambda$ in $\Sph_\infty^n$, let $\Sigma$ be any properly embedded minimal hypersurface such that   $\partial_\infty\Sigma=\Lambda$. Then $\Sigma$ is contained in  $\Cal(\Lambda)$. Moreover, if $\Lambda$ is not the boundary of a totally geodesic hyperplane, then $\Sigma$ is contained in the interior of $\Cal(\Lambda)$.
	\el
	\bp
	Let $P$ be any totally geodesic hyperplane disjoint from $\Lambda$. The signed distance function from $P$, defined in such a way that it goes to $-\infty$ as it approaches $\Lambda$, cannot have a positive maximum by a standard application of the geometric maximum principle (see also Corollary \ref{max principle} for a more general statement). This implies that $\Sigma$ is contained in the half-space bounded by $P$ whose closure contains $\Lambda$. Since $\Cal(\Lambda)$ is the intersection of all such half-spaces, this concludes the proof of the first assertion. 
	
	For the second assertion, suppose that $\Sigma$ contains a point $x$ in the boundary of $\Cal(\Lambda)$. Let $P$ be any support hyperplane for $\Cal(\Lambda)$ containing $x$. This means that $\Sigma$ is tangent to $P$ and, by the first assertion, it is contained in a half-space bounded by $P$. By the strong maximum principle, $\Sigma=P$, and therefore $\Lambda=\partial_\infty P$.
	\ep 
	
	\subsection{Stable and area-minimizing minimal hypersurfaces} A minimal hypersurface $\Sigma$ in a Riemannian manifold $M^{n+1}$ is {stable} if and only if, for every $u\in C_0^\infty(\Sigma)$,  
	$$\int_{\Sigma} uL_{\Sigma}(u)\mathrm{dVol}_\Sigma\geq 0~,$$
	where $L_\Sigma$ is the Jacobi operator:
	\begin{equation}\label{eq:jacobi}
		L_\Sigma(u)=-\Delta_\Sigma u-\left(\|A_\Sigma\|^2+\mathrm{Ric}_M(N_\Sigma,N_\Sigma)\right)u~.
	\end{equation}
	Equivalently, integrating by parts, $\Sigma$ is stable if and only if 
	\begin{equation}\label{eq:stability}
		\int_\Sigma\left(\|A_\Sigma\|^2+\mathrm{Ric}_M(N_\Sigma,N_\Sigma)\right)u^2\mathrm{dVol}_\Sigma\leq \int_\Sigma\|\nabla u\|^2\mathrm{dVol}_\Sigma
	\end{equation}
	for every $u\in C_0^\infty(\Sigma)$.
	
	A compact hypersurface $\Sigma_c$ with boundary is a \emph{least area hypersurface} if its area is less than or equal to that of any other compact hypersurface $\Sigma_c'$ such that $\partial\Sigma_c=\partial\Sigma_c'$. For the non-compact case, we say a hypersurface $\Sigma$ is area-minimizing if any compact codimension 0 submanifold with boundary of $\Sigma$ has least area in the sense above. 
	
	An area-minimizing hypersurface is stable. Indeed, the area-minimizing condition implies that $\Sigma$ is a critical point of the area functional among compactly supported variation (that is, $\Sigma$ is minimal) and the second derivative of the area is non-negative (which is equivalent to the definition of stability given above).

	When $n=2$, the following theorem provides an important a priori bound on the curvature of a stable minimal surface. We state the theorem here in its version for embedded minimal surfaces in $\HH$, but it holds more generally for immersed CMC surfaces in a complete Riemannian manifold of bounded sectional curvature. 
	
	\bt{{\cite[Main Theorem]{RST10}}}\label{thm:rosenberg}
	There exists a constant $c>0$ such that, for any stable embedded two-sided minimal surface $\Sigma$ in $\HH$, 
	$$\|A_\Sigma(p)\|\leq \frac{c}{\min\{d(p,\partial\Sigma),\pi/2\}}~.$$
	\et
	
	Applying Theorem \ref{thm:rosenberg} to a properly embedded minimal surface, which is automatically two-sided (see \cite{Sam69}), we immediately get the following corollary. See also \cite[Corollary 11]{Ros06}.
	
	\bcor \label{cor:rosenberg prop emb}
	There exists a constant $c'>0$ such that, for any stable properly embedded minimal surface $\Sigma$ in $\HH$, $\|A_\Sigma\|\leq c'$.
	\ecor

	\section{The Asymptotic Plateau Problem For Disks}\label{sec:APP}
	
	Anderson \cite[Theorem 4.1]{And83} stated the following important existence theorem for disk solutions to the asymptotic Plateau problem:
	\begin{theorem} \label{andersonthm}
		Let $\Lambda$ be a Jordan curve in $\mathbb S^2 = \partial_{\infty} \mathbb{H}^3$. There exists a (properly) embedded area-minimizing minimal disk asymptotic to $\Lambda$. 
	\end{theorem}
	However, there is a serious gap in the argument of the proof (see Remark 3.1 in \cite{Cos16}). The problem is that, without a uniform area upper bound, the cut-and-paste argument to produce competitors can increase the genus. So, while it is possible to solve the {\app} for $\Lambda$ by following his arguments, to control the topology of the {\ms}s one needs some additional arguments. We outline here a proof of Theorem \ref{andersonthm}, building on ideas of Coskunuzer, Gabai, and Soma. In fact our approach proves the existence of a properly embedded disk solution. The point that the embedding is {\it proper} was not addressed in \cite{And83}, but it will be essential in applications. This theorem will be then used, after appropriate modifications, in two key arguments for this paper, namely for the proofs of Theorem \ref{thm:char uniqueness} and Theorem \ref{thm:uncountably}.
	
	\subsection{Minimal laminations}
	
	We will need to use the notion of \emph{minimal lamination}. A minimal lamination $\mathcal L$ is a closed subset of $\HH$ which is the disjoint union of connected minimal surfaces (called \emph{leaves}), and which is covered by charts $(U\subset\HH,\Phi:U\to\R^3)$ in which $\Phi(\mathcal L\cap U)$ is represented as $\R^2\cap C$ for $C$ a closed subset of $\R$, in such a way that all leaves of $\mathcal L$ inside $U$ are mapped by $\Phi$ to horizontal slices of the form $\R^2\times\{\star\}$. A properly embedded (connected) minimal surface is an example of a minimal lamination, consisting of only one leaf.
	
	There is a notion of convergence of a sequence of minimal surfaces $\Sigma_n$ (or of minimal laminations, although we won't need this formulation) to a minimal lamination $\mathcal L$. This means that $\mathcal L$ is the set of accumulation points of sequences $x_n$ in $\HH$ where $x_n\in\Sigma_n$, and moreover, whenever one has a subsequence $x_{n_k}\to x\in\mathcal L$, then one can find smooth embeddings of disks in $\Sigma_{n_k}$  around $x_{n_k}$ (or, more in general, in the leaf through $x_{n_k}$) converging $C^\infty$ to a smooth minimal disk containing $x$. See \cite[Section 3]{Gab97}, \cite[Appendix B]{CMIV}, \cite[Section 3]{BT16} or \cite{Ba24} for more details.

	If $\Lambda$ is a smooth Jordan curve, in the absence of a uniform area bound, Gabai (\cite{Gab97}) constructed a minimal laminations whose asymptotic boundary is $\Lambda$, and he conjectured that one can find a properly embedded area-minimizing disk with any smooth asymptotic boundary $\Lambda$. This conjecture was proved by Soma (\cite{Som04}). This is our starting point. 
	
	\subsection{Monotone sequences}
	
	Inspired by \cite{Cos08,coskunuzer2019nestedsequences}, we prove Theorem \ref{andersonthm} by using (monotone) smooth approximations. Let us introduce some terminology. We say that a sequence of Jordan curves $\Lambda_n$ is \emph{monotone} if $\Lambda_n=\partial\Omega_n$, where $\Omega_n$ is a connected component of $\mathbb S^2\setminus\Lambda_n$, such that $\Omega_n\subseteq \Omega_{n+1}$. Similarly, we say that a sequence $\Sigma_n$ of properly embedded disks in $\HH$ is \emph{monotone} if $\Lambda_n=\partial\Delta_n$, where $\Delta_n$ is a connected component of $\HH\setminus\Sigma_n$, such that $\Delta_n\subseteq \Delta_{n+1}$. (Recall that, by the Jordan-Brouwer separation theorem, $\Sigma_n$ disconnects $\HH$ in two connected components.)
	
	Now, the fundamental ingredient will be the following convergence lemma, which itself can be of independent interest:
	
	\begin{lem}\label{lemma monotone convergence}
		Let $\Sigma_n$ be a monotone sequence of properly embedded minimal disks in $\HH$, such that $\Lambda_n=\partial_\infty\Sigma_n$ is a monotone sequence of Jordan curves converging $C^0$ to a Jordan curve $\Lambda$. Assume moreover that, for every compact subset $K\subset\HH$, there exists a constant $C=C(K)$ such that
		$$\sup_{\Sigma_n\cap K}\|A_{\Sigma_n}\|\leq C(K)~.$$
		Then $\Sigma_n$ converges $C^\infty$ to a properly embedded minimal disk $\Sigma$ in $\HH$ such that $\partial_\infty\Sigma=\Lambda$.
	\end{lem}
	
	\bp
	By \cite[Proposition B.1]{CMIV} (see also \cite[Lemma 3.3]{Gab97}), up to extracting a subsequence, the minimal disks $\Sigma_n$ converge  to a minimal lamination $\mathcal L$.
	
	Observe that the limiting lamination $\mathcal{L}$ might a priori be empty, but this cannot be the case here. Indeed, let $\gamma$ be a geodesic of $\HH$ whose endpoints at infinity are contained in different connected components of $\mathbb S^2_\infty\setminus\Lambda$.  If $n\geq n_0$, then one endpoint of $\gamma$ is contained in $\Omega_n$, while the other is always contained in the complement of $\Omega_n$. Since $\Sigma_n$ is properly embedded, by the Jordan-Brouwer separation theorem it disconnects $\HH$, which implies that $\Sigma_n$ must intersect $\gamma$. Let $x_n$ be the first intersection point on $\gamma$ when traveling from $\Omega_n$ to its complement. By monotonicity, $\{x_n\}_{n\in\N}$ is a monotone sequence on $\gamma$. 
	Moreover, by Lemma \ref{lemma interior convex hull}, every $\Sigma_n$ is contained in the convex hull $\mathcal C(\Lambda_n)$. Since $\Lambda_n\to\Lambda$, we have $\mathcal C(\Lambda_n)\to \mathcal C(\Lambda)$. Hence the sequence $\{x_n\}_{n\in\N}$ has an accumulation point on $\gamma\cap \mathcal C(\Lambda)$, which is a compact interval. Therefore $\mathcal L$ cannot be empty. 
	
	Let us also show that the asymptotic boundary of $\mathcal L$ equals $\Lambda$.  To see this, take a point $x\in \Lambda$ and let $\gamma_k$ be a sequence of geodesics in $\HH$ whose endpoints are contained in different connected components of $\mathbb{S}^2_{\infty} \setminus \Lambda$ and both converge to $x$. 
	By the same argument as above, for every $k$ there exists a point $y_k\in\mathcal L\cap \gamma_k$, obtained as the limit of a sequence of points in $\Sigma_n\cap\gamma_k$. Taking the limit in $k$, we obtain that $y_k$ converges to $x\in\Lambda$. Using the monotonic convergence, it is also easy to see that every point in the boundary at infinity of $\mathcal L$ must be in $\Lambda$, and thus the asymptotic boundary of $\mathcal L$ equals $\Lambda$.
	
	
	If $\Delta_n$ denotes the connected components of the complement of $\Sigma_n$ so that $\Delta_n\subseteq \Delta_{n+1}$ as in the definition of monotone sequence, we have $\mathcal L=\partial\Delta\subset\HH$, where $\Delta=\cup_{n=1}^\infty \Delta_n$. Indeed, one easily shows that every (subsequential or not) limit of any sequence $x_i\in\Sigma_i$ must be contained in $\partial\Delta$, and conversely every point of $\partial\Delta$ is the limit of such a sequence $x_i$. Incidentally, although not necessary for the proof, we observe that $\Delta$ is homeomorphic to an open ball by \cite{monotoneunion}. 
	
	Now, we claim that $\mathcal L$ is connected (i.e. it consists of only one leaf $\Sigma$) and properly embedded. Hence $\Sigma=\mathcal L$ is the desired solution of the asymptotic Plateau problem. 
	
	To see that $\mathcal L$ is connected, it suffices to show that the complement $F:=\HH\setminus\Delta$ is connected (see \cite{mathstack}). Clearly this is equivalent to showing that $\overline{F}=F\cup\partial_\infty F$ is connected, since $\partial_\infty F$ cannot disconnect. But $\overline F=\cap_{n=1}^\infty\overline{F_n}$, where $F_n:=\HH\setminus\Delta_n$, and the $\overline{F_n}$ are nested, compact and connected. This is sufficient to infer that $\overline F$ is connected. Indeed, if $\overline F$ were contained in the union of two disjoint open sets $U_1$ and $U_2$ of $\R^3$ (using the ball model), then each $\overline{F_n}$ would contain a point $y_n$ of $\partial U_1$, and by compactness one can extract a subsequence converging to $y\in\overline{F}\cap \partial U_1$, giving a contradiction. 
	
	By construction, $\mathcal L$ is closed and consists of only one leaf, and moreover has locally uniformly bounded curvature as a consequence of our hypothesis on the uniform bounds on $\|A_{\Sigma_n}\|$ over any compact subset $K$. Hence it is a properly embedded minimal surface (see \cite[Theorem 4]{MRclosure} applied to the case of a closed, connected surface with locally uniformly bounded curvature).

	Because $\mathcal{L}$ is a smooth limit of embedded disks, where the convergence is uniform on compact sets, $\mathcal{L}$ can be exhausted by  open precompact disks, and is therefore homeomorphic to a disk by \cite{monotoneunion}. This concludes the proof.  
	\ep
	
	\subsection{Conclusion of the argument}
	
	We are now ready to conclude a proof of the existence of a properly embbeded disk solution to the {\app}. 
	
	\bp[Proof of Theorem \ref{andersonthm}]
	Let $\Lambda$ be a Jordan curve on $\mathbb{S}^2_{\infty}$, and take a monotone sequence of smooth Jordan curves $\Lambda_n$ that $C^0$ converge to $\Lambda$. This is easily constructed, for example, by taking a biholomorphism between $\mathbb D$ and a connected component $\Omega$ of $\mathbb S^2_\infty\setminus\Lambda$, and considering the image of the ball of radius $1-1/n$ in $\mathbb D$. 
	
	For each $\Lambda_n$, we solve the asymptotic Plateau problem by \cite{Som04} and find area-minimizing properly embedded disks $\Sigma_n$ with boundary at infinity $\Lambda_n$. Since each $\Sigma_n$ is area-minimizing, it is stable, and therefore by Corollary \ref{cor:rosenberg prop emb} there exists $C$ such that $\|A_{\Sigma_n}\|\leq C$ for every $n$. Moreover, the $\Sigma_n$ form a monotone sequence as a consequence of the area-minimizing condition (see \cite[Lemma 3.2]{Cos08}). Hence we can apply Lemma \ref{lemma monotone convergence} and conclude. Note that the solution $\Sigma$ provided by Lemma \ref{lemma monotone convergence} is area-minimizing, since it is a uniform limit over compact sets of the area-minimizing surfaces $\Sigma_n$.
	\ep

	
	\section{Uniqueness conditions}\label{sec3}
	
	In this section, we prove Theorems \ref{thm:char uniqueness}, \ref{thm:char uniqueness invariant} and \ref{main1}.
	
	\subsection{Proof of Theorem \ref{thm:char uniqueness}}
	
	The key result to prove Theorem \ref{thm:char uniqueness} is the following:
	
	\bt\label{main3}
	Let $\Lambda$ be a {\Jc} on ${\Sph}^2_{\infty} = \partial \HH$. Suppose  that there exists a properly embedded {\ms} in $\HH$, which is not a stable minimal disk, spanning $\Lambda$. Then $\Lambda$ spans at least two disjoint  stable minimal disks in $\HH$.  
	\et
	
	Before proving Theorem \ref{main3}, let us explain how this implies Theorem \ref{thm:char uniqueness}. By Theorem \ref{andersonthm}, we know that $\Lambda$ spans a stable (actually, area-minimizing) properly embedded minimal disk. So uniqueness among minimal surfaces clearly implies uniqueness among stable minimal disks. For the converse direction, suppose by contradiction that $\Lambda$ spans a unique stable minimal disk, but it spans several minimal surfaces. Hence it spans a minimal surface $U$ which is not a stable minimal disk (that is, either it is not stable, or it is not a disk, or both). By Theorem \ref{main3}, $\Lambda$ spans two distinct stable minimal disks, and this gives a contradiction.

	To prove Theorem \ref{main3}, we first need the following result. 
	
	\bl \label{lemma Plateau}
	Let $\Lambda$ be a {\Jc} on ${\Sph}^2_{\infty} = \partial \HH$, let $U$ be a properly embedded {\ms} in $\HH$ spanning $\Lambda$, and let $\Omega$ be a connected component of $\HH\setminus U$. Then there exists a properly embedded stable minimal disk $\Sigma$ in $\Omega\cup U$ spanning $\Lambda$. 
	\el
	
	\bp
	We work as in the proof of Theorem \ref{andersonthm}.
	Let $\Lambda_n$ be a monotone sequence of smooth Jordan curves contained in $\partial_\infty\Omega$ and disjoint from $\Lambda$, that converge to $\Lambda$. 
	
	We claim that, for each $\Lambda_n$, we can solve the asymptotic Plateau problem in $\Omega$, that is, we can find a properly embedded minimal surface $\Sigma_n$ contained in $\Omega$, which is area-minimizing in $\Omega$, such that $\partial_\infty\Sigma_n=\Lambda_n$. Assuming the claim, the same argument as in \cite[Lemma 3.2]{Cos08} shows that  $\Sigma_n$ are a monotone sequence, and since each $\Sigma_n$ is area-minimizing in $\Omega$, it is stable,  so that we have $\|A_{\Sigma_n}\|\leq C$ for some $C>0$ and for all $n$. Hence we can apply Lemma \ref{lemma monotone convergence} and find a properly embedded stable minimal surface $\Sigma$, homeomorphic to a disk, which is the limit of the $\Sigma_n$. Since every $\Sigma_n$ is contained in $\Omega$, the limit $\Sigma$ is contained in its closure, namely in $\Omega\cup U$.
	
	It remains to show the claim. This essentially follows from the arguments in \cite{Som04}, that we now outline. To produce a solution of the asymptotic Plateau problem in $\HH$ for a smooth Jordan curve $\hat\Lambda$, Soma takes a sequence of simple closed curves $\gamma_k$ contained in $\HH$ and converging to $\hat\Lambda$. Using the results of \cite[Theorem 1]{MY19} or \cite[Theorem 6.3]{HS88}, one can find a solution of the finite Plateau problem, namely a least area disk with boundary $\gamma_k$. Moreover, Soma needs to choose this solution to be contained in the complement of a number of suitable  half-spaces in $\HH$, which is possible because the results of \cite{MY19} or \cite{HS88} hold in a sufficiently convex Riemannian manifold, that is, roughly speaking, a Riemannian manifold with piecewise smooth boundary and such that every smooth component of the boundary is mean convex (plus an additional extension assumption which is trivially satisfied here, see \cite[Section 1]{MY19} or \cite[Section 6]{HS88}). Then Soma shows that these solutions of the finite Plateau problem converge to a properly embedded area-minimizing disk $\hat\Sigma$ such that $\partial_\infty\hat\Sigma=\hat\Lambda$.
	
	In the present setting, we first observe that, when $\hat\Lambda=\Lambda_n$ the sequence $\gamma_k$, that converges to $\Lambda_n$, is eventually contained in $\Omega$. Second, we need to pick the solutions to the finite Plateau problem to be contained not only in a number of suitably chosen half-spaces, but also in $\Omega$. This is possible since the intersection of $\Omega$ and of some half-spaces is piecewise smooth up to a small perturbation, by genericity of the transverse intersection, and each smooth component of the boundary has vanishing mean curvature, so the definition of sufficiently convex is again satisfied. The arguments of Soma then produce, by taking limits of these finite solutions, a properly embedded area-minimizing disk $\Sigma_n$, which is a priori contained in $\Omega\cup U$, but is actually contained in $\Omega$ by the strong maximum principle. This concludes the argument.
	\ep

	We now prove Theorem \ref{main3}.
	
	\bp[Proof of Theorem \ref{main3}]
	Suppose $U$ is a properly embedded {\ms} in $\HH$ spanned by $\Lambda$, which is not a stable minimal disk. By the Jordan-Brouwer separation theorem, $U$ disconnects $\HH$. Let $\Omega^{\pm}$ be the two connected components of $\HH\backslash U$. 
	By Lemma \ref{lemma Plateau}, there exists a properly embedded stable minimal disk $\Sigma^\pm$ contained in $\Omega^\pm\cup U$ with $\partial_\infty\Sigma^\pm=\Lambda$. We claim that $\Sigma^\pm$ is disjoint from $U$. Indeed, $\Sigma^\pm$ is contained by construction in $\Omega^\pm\cup U$, so by the strong maximum principle, if $\Sigma^\pm$ and $U$ intersect, then they are equal.  
	Since $\Sigma^\pm$ is a stable minimal disk, and  $U$ is not, they must be disjoint.  This shows that $\Sigma^+\subset\Omega^+$ and $\Sigma^-\subset\Omega^-$, hence they are in particular disjoint.
	\ep

	\subsection{Small curvature conditions}\label{subsec:small curvatures}
	
	Next, we will show that a small curvature condition can imply uniqueness for solutions to the {\app}. This will prove Theorem \ref{main1}. 
	Let us first introduce (in any dimension) several conditions of having small curvature.

	\begin{definition}
		Let $\Sigma$ be an immersed hypersurface in $\mathbb H^{n+1}$, and let $\lambda_i$ denote its principal curvatures, for $i=1,\ldots,n$. We say that:
		\begin{itemize}
			\item $\Sigma$ has  \emph{weakly small curvature} if $|\lambda_i(x)|\leq 1$ for every $x$ in $\Sigma$ and every $i=1,\ldots,n$;
			\item $\Sigma$ has  \emph{small curvature} if $|\lambda_i(x)|< 1$ for every $x$ in $\Sigma$ and every $i=1,\ldots,n$;
			\item $\Sigma$ has  \emph{strongly small curvature} if there exists $\epsilon>0$ such that $|\lambda_i(x)|\leq 1-\epsilon$ for every $x$ in $\Sigma$ and every $i=1,\ldots,n$.
		\end{itemize}
	\end{definition}

	In this article, will a key result on hypersurfaces of weakly small curvatures, namely Theorem \ref{exponential diffeo} below. This result is well known under the small curvatures assumptions, see \cite{Eps84} and \cite[Section 4]{ES22}. In Appendix \ref{appendix} we will explain how the arguments adapt under the weakly small curvatures assumption.

	\bt \label{exponential diffeo}
	Let $\Sigma\subset\HHnp$ be an embedded hypersurface of weakly small  curvatures, and consider the map $F_\Sigma:\Sigma\times\R\to\HHnp$ defined by 
	$$F_\Sigma(x,t)=\exp_x(t N(x))~.$$
	\begin{enumerate}
		\item If $\Sigma$ is geodesically convex, then $F_\Sigma$ is a diffeomorphism onto its image, and for every $t\in\R$, $\Sigma_t:=F_\Sigma(\Sigma,t)$ is an embedded hypersurface of weakly small curvature. 
		\item If $\Sigma$ is complete, then $F_\Sigma$ is a diffeomorphism, and for every $t\in\R$, $\Sigma_t:=F_\Sigma(\Sigma,t)$ is a properly embedded hypersurface of weakly small curvature. 
	\end{enumerate}
	
	If moreover $\Sigma$ is minimal, then $H(\Sigma_t)$, the mean curvature of $\Sigma_t$, computed with respect to the unit normal vector field  pointing towards the direction where $t$ is increasing, has the same sign as $-t$ except possibly on a zero measure subset where $H(\Sigma_t)$ can vanish. 
	\et
	
	A consequence of Theorem \ref{exponential diffeo} is that, if $\Sigma$ is an embedded and geodesically convex hypersurface of weakly small  curvature,  the signed distance function $f_\Sigma:\mathcal N_\Sigma\to\R$ is smooth and has no critical points, where $\mathcal N_\Sigma:=F_\Sigma(\Sigma\times\R)$. Indeed, $f_\Sigma$ is simply the composition \begin{equation}\label{eq:signed distance}
		f_\Sigma=\pi\circ F_\Sigma^{-1}
	\end{equation}
	where $\pi:\Sigma\times\R\to\R$ is the projection to the second factor.

	Let us focus again on minimal hypersurfaces. The following statement is an immediate consequence of great importance for this work.
	
	\bcor \label{max principle}
	Let $\Sigma\subset\HHnp$ be an embedded geodesically convex minimal hypersurface of weakly small  curvature. For any embedded minimal hypersurface $\Sigma'\subset \mathcal N_\Sigma$, the restriction of  $d_{\Sigma}:=d(\cdot,\Sigma)$ to $\Sigma'$ has no positive local maximum. 
	\ecor
	\bp
	Assume by contradiction that $(d_{\Sigma})|_{\Sigma'}$ has a positive local maximum. Let $x_{\max}\in\Sigma'$ the maximum point and $t_{\max}=d_\Sigma(x_{\max})>0$ be the corresponding maximum value. Up to switching the sign of $f_\Sigma$, $d_{\Sigma}$ coincides with $f_\Sigma$ in a neighborhood of $x_{\max}$, where $f_\Sigma$ is the signed distance function as in \eqref{eq:signed distance}. Hence $(f_{\Sigma})|_{\Sigma'}$ also has a positive local maximum at $x_{\max}$. This implies that $\Sigma'$ is, in an embedded neighborhood of $x_{\max}$, tangent to $\Sigma_{t_{\max}}$ and contained in the region $F_\Sigma(\Sigma\times\{t\leq t_{\max}\})$, which is mean convex by the last part of Theorem \ref{exponential diffeo}. By the strong maximum principle for mean curvature, $\Sigma_{t_{\max}}$ coincides with $\Sigma'$, and hence has zero mean curvature, in a neighborhood of $x_{\max}$. But this contradicts the last claim of Theorem \ref{exponential diffeo}, namely that the subset of $\Sigma_{t_{\max}}$ where the mean curvature vanishes has zero measure.
	\ep

	Also, we mention the following result, by which we can restrict to properly \emph{embedded} hypersurfaces, when they have weakly small curvatures, as opposed to properly \emph{immersed}. See again Appendix \ref{appendix}.
	
	\bt \label{thm:weakly small is properly emb}
	Let $\Sigma\subset\mathbb H^{n+1}$ be a properly immersed hypersurface of weakly small  curvature. Then $\Sigma$ is properly embedded and diffeomorphic to $\R^n$. 
	\et

	Finally, it is immediate to check that minimal hypersurfaces of weakly small curvatures are stable.
	
	\bl\label{lemma weakly small is stable}
	Let  $\Sigma$ be a minimal hypersurface in $\HHnp$ of weakly small curvatures. Then $\Sigma$ is stable. 
	\el
	\bp
	If $\Sigma$ has weakly small curvatures, then $$\|A_\Sigma\|^2=\lambda_1^2+\ldots+\lambda_n^2\leq n=-\mathrm{Ric}_{\HHnp}(N,N)~.$$ Hence the left hand side of \eqref{eq:stability} is non-positive, and this concludes the proof.
	\ep

	\subsection{Proof of Theorem \ref{main1}} 
	

	In this subsection, we prove Theorem \ref{main1}. By Theorem \ref{thm:char uniqueness}, we know that it suffices to prove the uniqueness in the class of stable minimal disks. 
	
	The fundamental result is therefore the following theorem. 
	
	\bt\label{main2}
	Let $\Lambda$ be a {\Jc} on ${\Sph}^2_{\infty} = \partial \HH$ of finite width, and let $\Sigma$ be a properly embedded {\ms} in $\HH$ of weakly small curvature asymptotic to $\Lambda$. Then $\Sigma$ is the unique properly embedded stable {\ms} in $\HH$ asymptotic to $\Lambda$. 
	\et
	
	Combining Theorem \ref{main2} with Theorem \ref{thm:char uniqueness}, we immediately obtain the first part of Theorem \ref{main1}. The ``moreover'' part of Theorem \ref{main1} then follows from Theorem \ref{andersonthm}, which shows that every Jordan curve bounds at least one area-minimizing disk. 
	
	We begin by an easy lemma.
	
	\bl\label{lemma finite distance}
	Let $\Lambda$ be a {\Jc} on ${\Sph}^2_{\infty} = \partial \HH$ of width $w(\Lambda)<+\infty$, and let $\Sigma$ be any properly embedded surface contained in $\Cal(\Lambda)$ with $\partial_\infty\Sigma=\Lambda$ . Then for every $x\in \Cal(\Lambda)$, $d(x,\Sigma)\leq w(\Lambda)$. 
	\el
	\bp
	We may assume that $\Lambda$ is not the boundary of a totally geodesic hyperplane, as in that case the conclusion holds trivially. Let $x\in \Cal(\Lambda)$. By definition of the width, and since $\partial^\pm\Cal(\Lambda)$ is closed in $\HH$, there exist two geodesic segments joining $x$ to $y_{\pm}\in\Cal^\pm(\Lambda)$, each of length at most $w(\Lambda)$. Since $\Sigma$ is contained in the interior of $\Cal(\Lambda)$ by Lemma \ref{lemma interior convex hull}, and disconnects $\Cal(\Lambda)$ by the Jordan-Brouwer separation theorem, one of these two segments must meet $\Sigma$ at a point $y$, which is on the segment between $x$ and the other endpoint $y_\pm$. Hence
	$$d(x,\Sigma)\leq d(x,y)\leq d(x,y_\pm)\leq w(\Lambda)~.$$
	This concludes the proof.
	\ep
	
	\bp[Proof of Theorem \ref{main2}]
	Let $\Sigma'$ be a properly embedded stable minimal surface such that $\partial_\infty\Sigma'=\Lambda$. Consider the restriction to the surface $\Sigma'$ of the distance function $d_{\Sigma}:=d(\cdot,\Sigma)$. We want to show that $d_{\Sigma}$ vanishes identically on $\Sigma'$. This will imply $\Sigma'\subseteq \Sigma$ and thus, since $\Sigma'$ is also properly embedded, $\Sigma=\Sigma'$. 
	
	Suppose by contradiction that 
	\begin{equation}\label{eq sup contradiction}
		0<\sup_{\Sigma'}d_\Sigma\leq w(\Lambda)~.
	\end{equation}
	(The second inequality follows from  Lemma \ref{lemma finite distance}.) Let $\{x_n'\}_{n\in\mathbb N}$ be a maximizing sequence for $(d_\Sigma)|_{\Sigma'}$, namely


	\begin{equation}\label{eq sup distance}
		t_n:=d(x'_n,\Sigma) \ge \sup_{\Sigma'}d_\Sigma - \frac1n~.
	\end{equation}
	
	By Theorem \ref{exponential diffeo}, we can write $x_n'=\exp_{x_n}(t_n N)$, where $x_n\in\Sigma$ and $N(x_n)$ is a unit normal vector to $\Sigma$ at $x_n$.

	Now, fix a point $x_o\in\HH$. Let $\phi_n$ be isometries of $\HH$ such that $\phi_n(x_n) = x_o$. 
	We denote $\Sigma_n = \phi_n(\Sigma)$ and $\Sigma'_n = \phi_n(\Sigma')$. Observe that by construction $\Sigma_n$ contains $x_o$, while  $\Sigma_n'$ intersects the ball of radius $w(\Lambda)$ centered at $x_o$. 
	
	Recall that $\Sigma$ has weakly small curvatures, hence $\|A_{\Sigma}\|^2\leq 2$, and since this condition is preserved by isometries, $\|A_{\Sigma_n}\|^2\leq 2$.
	Since $\Sigma'$ is stable by assumption, and stability is preserved by isometries, $\Sigma'_n$ is stable for every $n$. Although in this case
	do not have explicit curvature bounds as for $\Sigma_n$, we use Corollary \ref{cor:rosenberg prop emb} to find that each $\|A_{\Sigma'_n}\|$ 
	is bounded by a universal constant. 
	
	Using elliptic regularity for the minimal submanifolds equation and the Ascoli-Arzelà Theorem for immersed submanifolds (see \cite{Smi07}), up to extracting a subsequence, we can assume that $\Sigma_n$ and $\Sigma_n'$ converge $C^\infty$ to immersed minimal submanifolds $\Sigma_\infty$ and $\Sigma'_\infty$ in a neighborhood of $x_o$ and of $x'_\infty=\lim_n x_n'$ respectively. Since the weakly small curvature condition is closed for the $C^\infty$ convergence, $\Sigma_\infty$ has weakly small curvature. Up to restricting further these neighborhoods if necessary, we can assume that $\Sigma_\infty$ is geodesically convex, so that by Theorem \ref{exponential diffeo} the exponential map is a diffeomorphism onto its image $\mathcal N_{\Sigma_\infty}$, and that $\Sigma_\infty'\subset \mathcal N_{\Sigma_\infty}$. Since every point of $\Sigma_n'$ can be written as $\exp_{\phi(y_n)}(s_n\phi_*(N(y_n)))$ for $s_n\leq \sup_{\Sigma'}d_\Sigma$, we have that the distance function $d_{\Sigma_\infty}=d(\cdot,\Sigma_\infty)$ restricted to $\Sigma_\infty'$ is bounded above by $D:=\sup_{\Sigma'}d_\Sigma$, and achieves the maximum value $D=d(x_o,x_\infty')$ at $x_\infty'$. This contradicts Corollary \ref{max principle} and thus concludes the proof.
	\ep




	\subsection{Proof of Theorem \ref{thm:char uniqueness invariant}} 
	
	We now provide a proof of Theorem \ref{thm:char uniqueness invariant}.  This is based on a combination of arguments in the spirit of the proof of Theorem \ref{thm:char uniqueness}  above, and the work of Guaraco, Lima and Vargas-Pallete \cite{GLP21}, which we summarize now.  Let $M$ be a quasi-Fuchsian manifold. In \cite[Theorem 1.5]{GLP21} it was proved that $M$ admits a foliation by closed surfaces that are either minimal, or have non-vanishing mean curvature (i.e., are either strictly mean-convex, or strictly mean-concave). The sign of the mean curvature switches at the minimal leaves.  
	This implies that if $M$ contains a unique closed stable minimal surface $S$, then 
	$M\setminus S$ has a foliation $\mathcal{F}$ by surfaces that are isotopic to $S$ and have mean-curvature vector everywhere pointing towards $S$.
	
	\begin{proof}[Proof of Theorem \ref{thm:char uniqueness invariant}]
		Suppose $\Gamma$ admits a unique invariant stable minimal disk $\Sigma$, which thus gives a unique closed stable minimal surface $\Sigma/\Gamma$ in $M:=\HH/\Gamma$. Let $F:M\to\R$ be a function whose level sets $F^{-1}(t)$ are the leaves of the foliation of $M$ provided by \cite{GLP21}, and such that $F^{-1}(0)=\Sigma/\Gamma$. The foliation $\mathcal{F}$ lifts to a foliation of $\HH$ by properly embedded disks $\widetilde\Sigma_t:=\widetilde F^{-1}(t)$ (where $\widetilde F$ is the lift of $F$ to $\HH$, and therefore $\widetilde\Sigma_0=\Sigma$) with asymptotic boundary the limit set $\Lambda$ of $\Gamma$.

		
		By Theorem \ref{thm:char uniqueness}, it suffices to show that, if   $\Sigma'$ is a properly embedded stable minimal disk  in $\HH$ with $\partial_\infty\Sigma'=\Lambda$, then $\Sigma'=\Sigma$. Consider the function $f:=\widetilde F|_{\Sigma'}:\Sigma'\to\R$. We claim that $f\equiv 0$, which will imply $\Sigma'=\Sigma$ since both are properly embedded. First, observe that $f$ is bounded. Indeed, since $\Gamma$ is quasi-Fuchsian, the convex core $\Cal(\Lambda)/\Gamma$ of $M$ is compact, and $F$ is thus bounded on the convex core. Lifting to $\HH$, this shows that $\widetilde F|_{\Cal(\Lambda)}$ is bounded, and so too is $f$ by Lemma \ref{lemma interior convex hull}.
		
		Now, the proof is similar to the proof of Theorem \ref{main2}, using the foliation lifted from $\mathcal F$ instead of the equidistant foliation. Suppose by contradiction that $\sup f>0$, the other case being analogous. Let $x_n'$ be a  sequence such that $f(x_n')\to \sup f$, let $\phi_n$ be a sequence of isometries sending $x_n'$ to a fixed point $x_o$, and let $\Sigma_n':=\phi_n(\Sigma')$ and $\Sigma_n:=\phi_n(\widetilde\Sigma_{\sup f})$. 
		
		Observe first that the width $\Cal(\Lambda)$ is finite because $\Lambda$ is a quasicircle, and, by cocompactness of the action of $\Gamma$ on $\widetilde\Sigma_{\sup f}$, the distance of $\widetilde\Sigma_{\sup f}$ from $\Cal(\Lambda)$ is bounded by some constant $c$. By construction $\Sigma_n'$ contains $x_o$, and it follows that $\Sigma_n$ intersects every ball centered at $x_o$ of radius $w(\Lambda)$.

		Second, by smoothness of $\widetilde\Sigma_{\sup f}$ and cocompactness again, the norm of the second fundamental form of $\widetilde\Sigma_{\sup f}$ and of all its covariant derivative are  bounded. Since $\phi_n$ is an isometry, the norm of the second fundamental form of $\Sigma_n$ and of all its covariant derivative are bounded too. The same holds for $\Sigma_n'$ by Corollary \ref{cor:rosenberg prop emb} and standard elliptic regularity theory.
		
		By the Ascoli-Arzelà Theorem (see \cite{Smi07}), up to extracting a subsequence, we have that $\{\Sigma_n\}_{n\in\mathbb N}$ and $\{\Sigma_n'\}_{n\in\mathbb N}$ converge smoothly as immersed surfaces to, say, $\Sigma_\infty$ and $\Sigma_\infty'$ in small neighborhoods. As we are assuming $\sup f>0$, and by cocompactness, $\widetilde\Sigma_{\sup f}$ has uniformly negative mean curvature with respect to the unit normal pointing towards $\widetilde F=+\infty$, hence the same holds for the limit  $\Sigma_\infty$.

		It follows from the construction that $x_o\in\Sigma'_\infty$. By cocompactness, the distance between different leaves $\widetilde \Sigma_{\sup f}$ and $\widetilde \Sigma_{t}$ goes to zero as $t\to\sup f$. Using this fact, one can show similarly to the proof of Theorem \ref{main2} that $x_o\in\Sigma_\infty$. Moreover, since $\Sigma'$ is contained in the region $\{x\,|\,\widetilde F(x)\leq \sup f\}$ whose boundary is $\widetilde\Sigma_{\sup f}$, after taking limits $\Sigma'_\infty$ is contained in the side of $\Sigma_\infty$ corresponding to decreasing values of $\widetilde F$. This contradicts the geometric maximum principle at the tangency point $x_o$, and therefore shows that $\sup f\leq 0$. Repeating the argument for the infimum of $f$ instead of the supremum shows that $f\equiv 0$ and concludes the proof.
	\end{proof}

	\section{Higher dimension and codimension}\label{sec:higher}
	In this section, we discuss some generalization of previous results to higher dimensions and higher codimensions.

	\subsection{Hypersurfaces in higher dimensions}
	
	We provide here a higher dimensional generalization of Corollary \ref{cor strongly small}.
	
	\bt\label{thm strongly small higher dim}
	Let $\Lambda$ be a topologically embedded $(n-1)$-sphere on ${\Sph}^n_{\infty} = \partial_{\infty} \HHnp$, and let $\Sigma$ be a properly embedded minimal hypersurface in $\HHnp$ of 
	{{strongly}} {\spc} asymptotic to $\Lambda$. Then $\Sigma$ is the unique properly embedded  minimal hypersurface in $\HHnp$ asymptotic to 
	$\Lambda$. Moreover, $\Sigma$ is area-minimizing.
	\et
		The argument is similar to the proof of Theorem \ref{main2}. However, in that case it was sufficient to work with a \emph{stable} minimal surface, by Theorem \ref{thm:char uniqueness}. The difference here is that instead of passing to a limit of minimal submanifolds (which we will be unable to do if they are unstable), we pass to a limit of metric balls in the ambient space.  We are then able to get a contradiction by applying the following theorem of White \cite{Whi10} --- which we state here for minimal hypersurfaces, but holds more generally in the setting of minimal varifolds.
		
		\begin{theorem}[{\cite[Theorem 3]{Whi10}}]
			\label{white}
			Let $M_n$ be a sequence of Riemannian manifolds with boundary, converging smoothly to a limit Riemannian manifold with boundary $M_\infty$. Suppose $\Sigma_n$ is a minimal hypersurface of $M_n$. Then no point of $\partial M_\infty$ where the boundary is strictly mean convex can be a limit of a sequence $x_n$ such that $x_n\in \Sigma_n$.
		\end{theorem}
		
		In order to apply Theorem \ref{white}, we need to work on a sequence $M_n$ with \emph{uniformly} strictly mean convex boundary (where the boundary is equidistant from the minimal hypersurface), and for this reason we need to make the assumption of \emph{strongly} small curvature. This assumption also ensures finite width (see Lemma \ref{lemma finite width}), which is another essential ingredient.

		\begin{proof}[Proof of Theorem \ref{thm strongly small higher dim}] 
			
			As in the proof of Theorem \ref{main2}, let $\Sigma'$ be a properly embedded  minimal hypersurface such that $\partial_\infty\Sigma'=\Lambda$, and consider the restriction to the surface $\Sigma'$ of the distance function $d_{\Sigma}:=d(\cdot,\Sigma)$. 
			Suppose by contradiction that 
			$
			0<\sup_{\Sigma'}d_\Sigma\leq w(\Lambda)$, where the upper bound 
			follows from  Lemma \ref{lemma finite distance}, that works identically in any dimension.
			
			Let $\{x_n'\}_{n\in\mathbb N}$ be a maximizing sequence for $(d_\Sigma)|_{\Sigma'}$
			By Theorem \ref{exponential diffeo} we can represent  $x_n'=\exp_{x_n}(t_n N)$, where $x_n\in\Sigma$ and $N(x_n)$ is a unit normal vector to $\Sigma$ at $x_n$. Observe that $t_n=d(x_n',\Sigma)\to D:=\sup_{\Sigma'}d_\Sigma$.

			Now, fix a point $x_o\in\HH$. Let $\phi_n$ be isometries of $\HH$ such that $\phi_n(x_n) = x_o$. 
			We denote $\Sigma_n = \phi_n(\Sigma)$ and $\Sigma'_n = \phi_n(\Sigma')$. Since $\Sigma$ has weakly small curvatures, $\|A_{\Sigma}\|^2\leq n$, and thus $\|A_{\Sigma_n}\|^2\leq n$. Using elliptic regularity for the minimal submanifolds equation and the Ascoli-Arzelà Theorem for immersed submanifolds (see \cite{Smi07}), up to extracting a subsequence, we can assume that $\Sigma_n$ converges $C^\infty$ in a ball $B(x_o,\delta)$ of radius $\delta$ centered at $x_o$ to an immersed (actually, embedded, up to taking a smaller $\delta$) minimal submanifold $\Sigma_\infty$. By the strongly small curvature assumption, the principal curvatures of $\Sigma_n$
			are in by $[-1+\epsilon,1-\epsilon]$ for some $\epsilon$, and thus the same holds for $\Sigma_\infty$.
			
			Up to restricting further $\delta$, we can assume that $\Sigma_\infty$ contains a geodesically convex neighborhood, and therefore by Theorem \ref{exponential diffeo} the exponential map is a diffeomorphism onto its image. 
			This implies that the $D$-equidistant hypersurface $\Sigma_\infty^D$ from $\Sigma_\infty$, that contains $y_o:=\exp_{x_o}(DN_{\infty}(x_o))$ and is the image of the map $\iota_{D}$ (see Lemma \ref{lemma principal curvatures}), is embedded and its mean curvature is bounded below by a positive constant (that depends on $D$ and $\epsilon$). 
			
			We are now ready to apply Theorem \ref{white}. Let $\Sigma_n^{t_n}$ be the $t_n$-equidistant surface from $\Sigma_n$ which is converging to $\Sigma_\infty^D$. Define $M_n$ to be the intersection between the open ball $B(y_o,\delta')$ and the closure of the connected component of $\HHnp\setminus\Sigma_n^{t_n}$ which is mean convex, where $\delta'$ is chosen sufficiently small so that $\Sigma_n^{t_n}$ converges smootly to $\Sigma_\infty^D$ inside $B(y_o,\delta')$ (the smooth convergence is clear, since we are assuming that $\Sigma_n$ converges to $\Sigma_\infty$ inside $B(x_o,\delta)$). Now, $M_n$ is an (incomplete) Riemannian manifold with boundary, as in Theorem \ref{white}. Its limit $M_\infty$ is the intersection of $B(y_o,\delta')$ with the closure of the connected component of $\HHnp\setminus\Sigma_\infty$ which is strictly mean convex. But $\Sigma_n':=\phi_n(\Sigma')$ is a sequence of minimal hypersurfaces in $M_n$ going through $\phi_n(x_n')$, and by construction $\phi_n(x_n')\to y_o\in\partial M_\infty$, contradicting Theorem \ref{white}.
			
			For the ``moreover'' part, Anderson \cite{And83} proved the existence of an area-minimizing integral current with asymptotic boundary $\Lambda$, which, however, is a smooth hypersurface only in the complement of a singular set of Hausdorff dimension at most $n-7$. Nevertheless, calling this integral current $\Sigma'$, the above argument based on Theorem \ref{white} applies to such a $\Sigma'$ as well (see \cite{Whi10}), showing that $\Sigma'=\Sigma$ is, in particular, a smooth hypersurace. Hence $\Sigma$ is area-minimizing.
		\end{proof}

		\subsection{Higher codimension submanifolds}
		
		We briefly outline some generalizations to higher codimension. By following Uhlenbeck's approach \cite{Uhl83},  Jiang \cite[Lemma 2.2]{j21} defined a notion of almost Fuchsian for quotients $M$ of $\mathbb{H}^{n+1}$ by discrete faithful representations of surface groups in $\text{Isom}(\mathbb{H}^{n+1})$. 
		Such a quotient $M$ is called almost-Fuchsian if it admits a (two-dimensional) minimal surface $\Sigma$ with principal curvatures smaller than 1 in magnitude.  Jiang showed that an almost-Fuchsian $M$ with $\Sigma $ removed is foliated by equidistant surfaces to $\Sigma$ that are mean 2-convex (the sum of the smallest two principal curvatures is positive), and thus by the maximum principle that $\Sigma$ is the unique closed minimal surface in $M$.  Bronstein recently showed that in four dimensions $M$ could in fact be homeomorphic to a nontrivial disk bundle over a surface \cite{b23}.  The limit set of such an $M$ would then be a self-similar everywhere-non-tame knot in   $\partial_\infty \mathbb{H}^4 \cong \mathbb S^3$ (see \cite{glt88}.) 
		
		For any almost-Fuchsian $M$, the foliation by mean 2-convex equidistant surfaces to $\Sigma$ lifts to a foliation of the complement of a lift $\widetilde{\Sigma}$ of $\Sigma$ to $\mathbb{H}^{n+1}$.  We thus expect that the same argument by contradiction using \cite[Theorem 3]{Whi10} (which applies in any codimension) can be used,  to show that the limit set of $M$ bounds a unique minimal disk, namely $\widetilde{\Sigma}$.  Jiang's work is for two dimension surfaces in arbitrary codimension, but we expect that similar arguments would work for $k$-dimensional minimal surfaces with small principal curvatures in any dimension. In that case, the same arguments would give a uniqueness result for asymptotic Plateau problems in that setting as well.

		Finally here we mention that Fine has recently developed a theory for counting minimal surfaces in $\mathbb{H}^4$ with asymptotic boundary a given knot \cite{f21}.  In view of \cite{b23} and the previous discussion, one should  be able to give examples of (wild)  knots in $\mathbb S^3$ that bound a unique minimal surface in $\mathbb{H}^4$.

		
		\section{Uncountably many minimal disks}\label{sec:uncountably}
		
		The aim of this section is to prove Theorem \ref{thm:uncountably}.
		\subsection{Constructing the Jordan curve}\label{subsec:construction jordan curve}
		
		We begin by constructing the Jordan curve $\Lambda$ in $\Sph^2_\infty=\partial_\infty\HH$, that will be the asymptotic boundary of uncountably many stable minimal disks. For this purpose, we will construct an explicit continuous injective map 
		$$f:\R\mathrm P^1=\R\cup\{\infty\}\to \C\mathrm P^1=\C\cup\{\infty\}~.$$
		
		We begin by fixing some notation. Let $\epsilon,\delta>0$ be some constants  which will be fixed later. For the moment, we only need to assume that $\epsilon<1/3$, so that the circle of radius $2^n(1+\epsilon)$ is smaller than the circle of radius $2^{n+1}(1-\epsilon)$, and $\delta<\pi/8$. Let 
		\begin{equation*}
			\centering
			\begin{aligned}
				p^n_{+,+}:=&2^n\left(1-{\epsilon}\right)e^{ i\left(\frac{\pi}{2}-\delta\right)} \\
				q^n_{+,+}:=&2^n\left(1+{\epsilon}\right)e^{ i\left(\frac{\pi}{2}-\delta\right)} \\
				p^n_{-,+}:=&2^n\left(1-{\epsilon}\right)e^{ i\left(\frac{\pi}{2}+\delta\right)}    \\
				q^n_{-,+}:=&2^n\left(1+{\epsilon}\right)e^{ i\left(\frac{\pi}{2}+\delta\right)}    
			\end{aligned}
			\qquad
			\begin{aligned}
				p^n_{+,-}:=&2^n\left(1-{\epsilon}\right)e^{ i\left(-\frac{\pi}{2}+\delta\right)}  \\
				q^n_{+,-}:=&2^n\left(1+{\epsilon}\right)e^{ i\left(-\frac{\pi}{2}+\delta\right)}  \\
				p^n_{-,-}:=&2^n\left(1-{\epsilon}\right)e^{ i\left(-\frac{\pi}{2}-\delta\right)}  \\
				q^n_{-,-}:=&2^n\left(1+{\epsilon}\right)e^{ i\left(-\frac{\pi}{2}-\delta\right)}  \\
			\end{aligned}
		\end{equation*}
		
		Observe that the points $p^n_{\sigma_1,\sigma_2}$ and $q^n_{\sigma_1,\sigma_2}$ lie in the quadrant determined by the two signs $\sigma_1$ and $\sigma_2$, and that the Euclidean distance between $p^n_{\sigma_1,\sigma_2}$ and $q^n_{\sigma_1,\sigma_2}$ is $2^{n+1}\epsilon$. Moreover, denoting 
		$\alpha(z)=2z$, we see that $p^n_{\sigma_1,\sigma_2}=2^np^0_{\sigma_1,\sigma_2}$
		and the same for $q^n_{\sigma_1,\sigma_2}$.
		
		Now, we define the map $f$, which will be injective and piecewise smooth when restricted to $\R\setminus\{0\}$, with junctions exactly at the points $p^n_{\sigma_1,\sigma_2}$ and $q^n_{\sigma_1,\sigma_2}$.  The explicit formula for $f$ is the following:
		
		$$f(t)=\begin{cases} 
			0 &\textrm{ if } t=0 \\
			2^n\left(1-{\epsilon}\right)e^{ i\left((4\pi-8\delta)t-\frac{9\pi}{2}+9\delta\right)} &\textrm{ if } t\in \left[2^n,2^n\left(1+\frac{1}{4}\right)\right] \\ 
			2^n\left(1+\epsilon(-8t+7)/3\right)e^{ i\left(\frac{\pi}{2}-\delta\right)} &\textrm{ if } t\in \left[2^n\left(1+\frac{1}{4}\right),2^n\left(1+\frac{1}{2}\right)\right] \\ 
			2^n\left(1+{\epsilon}\right)e^{ i\left((-4\pi+8\delta)t+\frac{5\pi}{2}-5\delta\right)} &\textrm{ if } t\in \left[2^n\left(1+\frac{1}{2}\right),2^n\left(1+\frac{3}{4}\right)\right] \\ 
			2^n\left((-4/3+4\epsilon)t+(10/3-6\epsilon)\right)e^{ i\left(-\frac{\pi}{2}+\delta\right)} &\textrm{ if } t\in \left[2^n\left(1+\frac{3}{4}\right),2^{n+1}\right] \\ 
			\infty &\textrm{ if } t=\infty \\
			-2^n\left(1-{\epsilon}\right)e^{ i\left((4\pi-8\delta)t+\frac{9\pi}{2}-9\delta\right)} &\textrm{ if } t\in \left[-2^n\left(1+\frac{1}{4}\right),-2^n\right] \\ 
			-2^n\left(1+\epsilon(8t+7)/3\right)e^{ i\left(-\frac{\pi}{2}+\delta\right)} &\textrm{ if } t\in \left[-2^n\left(1+\frac{1}{2}\right),-2^n\left(1+\frac{1}{4}\right)\right] \\ 
			-2^n\left(1+{\epsilon}\right)e^{ i\left((-4\pi+8\delta)t-\frac{5\pi}{2}+5\delta\right)} &\textrm{ if } t\in \left[-2^n\left(1+\frac{3}{4}\right),-2^n\left(1+\frac{1}{2}\right))\right] \\ 
			-2^n\left((4/3-4\epsilon)t+(10/3-6\epsilon)\right)e^{ i\left(\frac{\pi}{2}-\delta\right)} &\textrm{ if } t\in \left[-2^{n+1},-2^n\left(1+\frac{3}{4}\right)\right] \\ 
		\end{cases}$$
		
		\begin{figure}[htb]
			\centering
			\includegraphics[width=.8\textwidth]{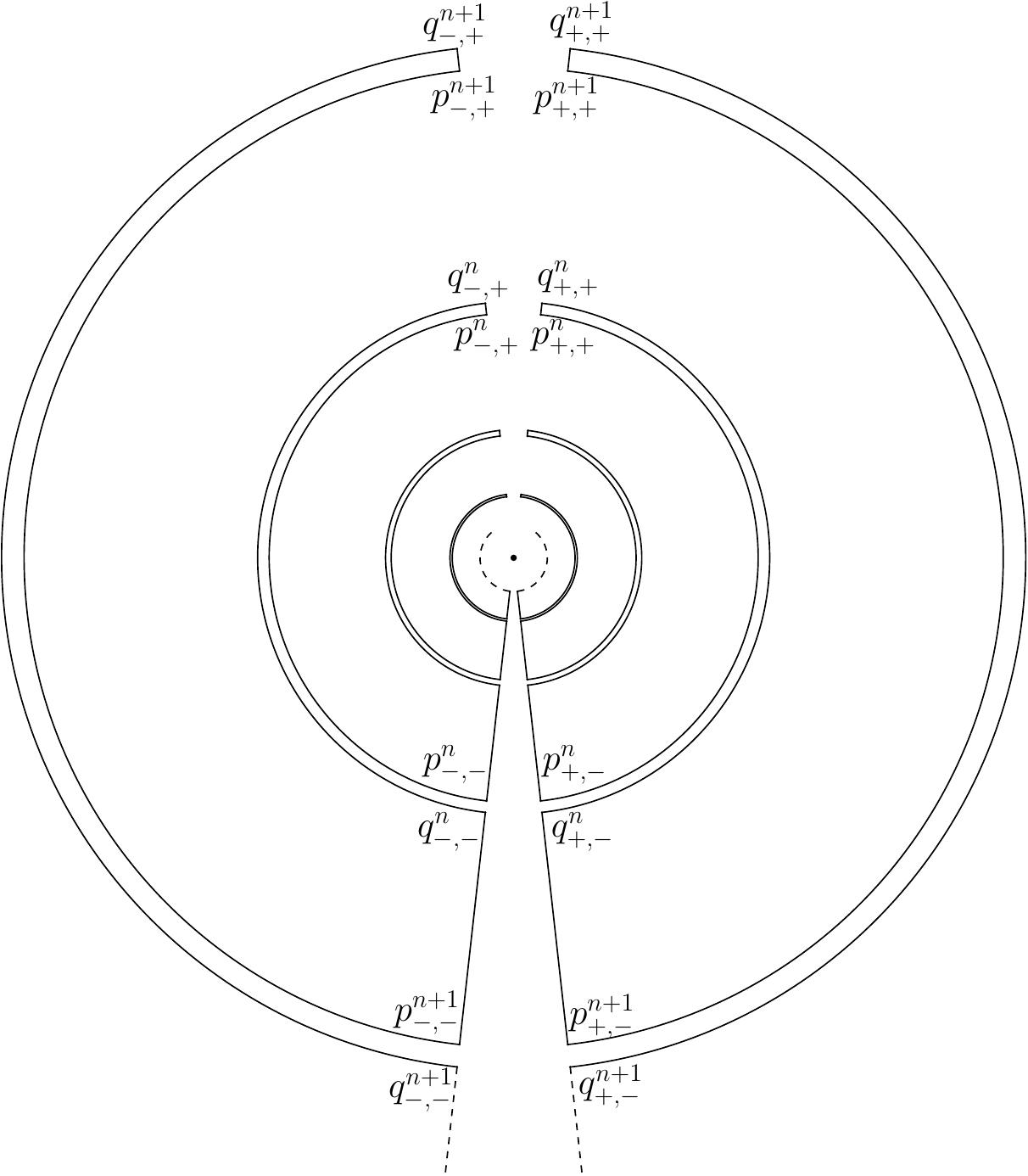} 
			\caption{\small The construction of the Jordan curve $\Lambda$.}\label{fig:Jordancurve}
		\end{figure}

		Let us explain this definition (see also Figure \ref{fig:Jordancurve}). First, observe that the image of $(0,+\infty)$ lies in the half-space $\{\Re(z)>0\}$, while the image of $(-\infty,0)$ lies in the half-space $\{\Re(z)<0\}$. The intervals $[2^n,2^n(1+1/4)]$ are mapped to the 
		arcs of circles contained in $\{\Re(z)>0\}$ connecting $p^n_{+,-}$ to $p^n_{+,+}$, and similarly the intervals $[2^n(1+1/2),2^n(1+3/4)]$ are mapped to the 
		arcs of circles contained in $\{\Re(z)>0\}$ connecting $q^n_{+,+}$ to $q^n_{+,-}$, parameterized proportionally to arclenght. 
		The intervals $[2^n(1+1/4),2^n(1+1/2)]$ are mapped to straight lines connecting $p^n_{+,+}$ to $q^n_{+,+}$. The intervals $[2^n(1+1/4),2^n(1+1/2)]$ are mapped to straight lines connecting $q^n_{+,-}$ to $p^{n+1}_{+,-}$. For negative values of $t$, the situation is absolutely analogous up to composing with a reflection in the imaginary axis. Namely, $f$ satisfies the symmetry: 
		\begin{equation}
			\label{symmetry f}
			f(-t)=-\overline{f(t)}~.
		\end{equation}
		Also, $f$ satisfies the equivariance property 
		\begin{equation}
			\label{equivariance f}
			f\circ \alpha^n=\alpha^n\circ f~,
		\end{equation}
		where we recall that $\alpha(z)=2z$.

		\bl\label{lemma quasicircle}
		The image of $f$ is a quasicircle in $\C\mathrm P^1$.
		\el
		\bp
		We show that $f$ is continuous and injective, which shows that the image of $f$ is a Jordan curve. One immediately checks from the formula that $f$ is well-defined and continuous (actually, piecewise smooth) when restricted to $\R\setminus\{0\}$. Observe moreover that, if $|t|\in [2^n,2^{n+1}]$, then $|f(t)|\in [2^n(1-\epsilon),2^{n+1}(1-\epsilon)]$. This immediately implies that $f$ is continuous at $t=0$ and $t=\infty$. 
		
		To check that $f$ is injective, we have already observed that $f$ maps $(0,+\infty)$ to the open right half-space and $(-\infty,0)$ to the open left half-space. Hence, together with \eqref{symmetry f}, it suffices to check that $f$ is injective on $(0,+\infty)$. By the condition on $|f(t)|$ in the previous paragraph and \eqref{equivariance f} it is actually sufficient to check that $f$ is injective on $[1,2]$, which is an immediate consequence of the construction.
		
		Let us denote by $\Lambda$ the image of $f$. We now show that $\Lambda$ is a quasicircle, by using the following characterization of quasicircles given in \cite[Chapter IV, Theorem 5]{Ahl66}. Let $\Lambda$ be a Jordan curve in $\C\mathrm P^1=\C\cup\{\infty\}$ going through $\infty$. Then $\Lambda$ is a quasicircle if and only if there exists a constant $C>0$ such that, for every $z_1,z_2\in\Lambda\cap \C$ and every $z_3$ in the arc of  $\Lambda\cap \C$ with endpoints $z_1,z_2$, 
		\begin{equation}\label{eq ahlfors}
			|z_3-z_1|\leq C|z_2-z_1|~.
		\end{equation}
		To check that \eqref{eq ahlfors} holds for the image of $f$, first observe that the condition \eqref{eq ahlfors} is scale-invariant. Hence, since $\Lambda$ is invariant by the transformations $\alpha^n(z)=2^nz$, we can assume that $\max\{|z_1|,|z_2|\}\in [1/2,1]$. By definition of $f$, the arc of  $\Lambda\cap \C$ with endpoints $z_1,z_2$ is then contained in the Euclidean ball $B(0,1)$. It follows that it is enough to check  \eqref{eq ahlfors} for $|z_2-z_1|\leq \epsilon_0$, for any fixed $\epsilon_0>0$. Indeed, if $|z_2-z_1|\geq \epsilon_0$, since by the triangle inequality $|z_3-z_1|\leq 2$, we have $|z_3-z_1|\leq (2/\epsilon_0)|z_2-z_1|$. In conclusion, it is sufficient to check that \eqref{eq ahlfors} holds when either $z_1$ or $z_2$ has modulus in $[1/2,1]$ and $|z_2-z_1|\leq \epsilon_0$. This is immediately verified since $\Lambda\cap \{1/2-\epsilon_0\leq |z|\leq 1+\epsilon_0\}$ is a union of finitely many straight line segments or arcs of circles meeting orthogonally at the adjacency points.
		\ep
		
		\begin{rem}
			We remark that $\Lambda$ is not a symmetric quasicircle. Indeed, as a consequence of the invariance under the scaling $\alpha(z)=2z$, $\Lambda$ cannot satisfy the so-called strong reverse triangle property in a neighborhood of $0$, see \cite[Section 6]{GS92}.
		\end{rem}
		
\begin{figure}[htb]
			\centering
			\includegraphics[width=.3\textwidth]{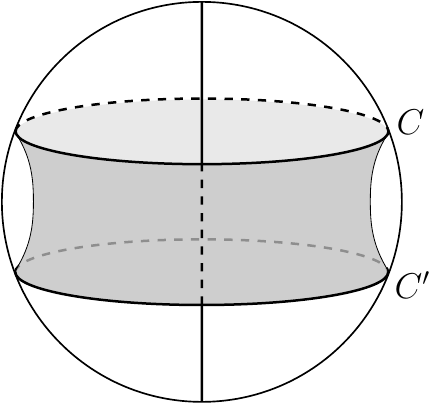} 
			\caption{\small 
				A minimal catenoid in $\HH$ and its axis of revolution.}\label{fig:catenoid}
		\end{figure}

		Now, in order to fix the suitable values of $\epsilon$ and $\delta$, we need to introduce minimal catenoids. A \emph{minimal catenoid} is a minimal surface in $\HH$, homeomorphic to an open annulus, whose asymptotic boundary is the  union of two disjoint circles $C$ and $C'$ in $\partial_\infty\HH=\C\mathrm P^1$, and which is obtained as a surface of revolution with respect to the unique geodesic meeting orthogonally both totally geodesic planes $P$ and $P'$ with $\partial_\infty P=C$ and  $\partial_\infty P'=C'$. See Figure \ref{fig:catenoid}.

		We will use the following existence result:
		\bt[\cite{MG87}]\label{existence min catenoids}
		There exists a constant $d>0$ such that, if $C$ and $C'$ are disjoint circles in $\partial_\infty\HH=\C\mathrm P^1$ and the distance between the totally geodesic planes $P$ and $P'$ that they span is less than $d$, then there exists a minimal catenoid with asymptotic boundary equal to $C\cup C'$.
		\et
		
		Now, let us define the following circles in $\C\mathrm P^1=\C\cup\{\infty\}$:
		
		\begin{align*}
			C_{n,0}&=\left\{z\,:\,\left| z-2^n\left(1-{2\epsilon}-\frac{1+\epsilon}{100}\right)\right|=\frac{2^n}{100}\right\} \\
			C'_{n,0}&=\left\{z\,:\,\left| z-2^n\left(1+{2\epsilon}+\frac{1+\epsilon}{100}\right)\right|=\frac{2^n}{100}\right\} \\
			C_{n,1}&=\{z\,:\,|z-2^n e^{i\left(\frac{\pi}{2}-\delta-(1+\delta)\epsilon\right)}|={\epsilon}\} \\
			C'_{n,1}&=\{z\,:\,|z-2^n e^{i\left(\frac{\pi}{2}+\delta+(1+\delta)\epsilon\right)}|=\epsilon\} \\
		\end{align*}

\begin{figure}[htb]
			\centering
			\includegraphics[width=.5\textwidth]{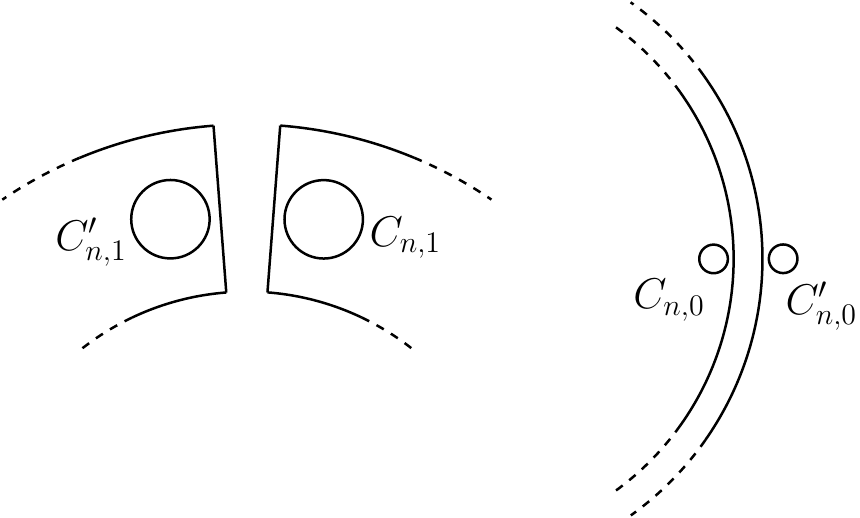} 
			\caption{\small The circles $C_{n,0},C_{n,0}',C_{n,1}$ and $C_{n,1}'$.}\label{fig:circles}
		\end{figure}

		See Figure \ref{fig:circles}. Observe that
		\begin{equation}
			\label{equivariance circles}
			C_{n,0}=\alpha^nC_{0,0}\qquad C'_{n,0}=\alpha^nC'_{0,0}\qquad C_{n,1}=\alpha^nC_{0,1}\qquad C'_{n,1}=\alpha^nC'_{0,1}~.
		\end{equation}
		and that $C_{n,0}$, $C_{n,0}'$, $C_{n,1}$ and $C_{n,1}'$ are all in the complement of the image of $f$.
		
		\bl
		There exists constants $\epsilon,\delta>0$ such that for all $n\in\mathbb Z$ and all $j\in\{0,1\}$, $C_{n,j}\cup C_{n,j}'$ is the asymptotic boundary of a minimal catenoid.
		\el
		\bp
		By \eqref{equivariance circles} and the fact that $\alpha$ is the boundary value of an isometry of $\HH$ in the upper half-space model, it suffices to prove the statement for $n=0$. We shall first choose $\epsilon$ in such a way that $C_{0,0}\cup C_{0,0}'$ bounds a minimal catenoid. For this, it suffices to observe that $C_{0,0}$ and $C_{0,0}'$ are circles of radii $1/100$, and the Euclidean distance between their centers is $2(2\epsilon+(1+\epsilon)/100)=2/100+402\epsilon/100$. As $\epsilon\to 0$, the circles $C_{0,0}$ and $C_{0,0}'$ tend to a pair of tangent circles, each with radius $1/100$. Since the hyperbolic distance varies continuously, and the hyperbolic distance between two totally geodesic planes tangent at infinity equals zero, the distance between the totally geodesic planes bounded by $C_{0,0}$ and $C_{0,0}'$ tends to zero. Hence by Theorem \ref{existence min catenoids}, when $\epsilon$ is sufficiently small, there exists a minimal catenoid with asymptotic boundary $C_{0,0}\cup C_{0,0}'$.
		
		We shall now keep such $\epsilon$ fixed, and choose $\delta$ in such a way that $C_{0,1}\cup C_{0,1}'$ bounds a minimal catenoid. The reasoning is completely analogous to the previous case. Indeed, the circles $C_{0,1}$ and $C_{0,1}'$ have radius $\epsilon$ and the distance between their centers is bounded above by $2\epsilon+2\delta(1+\epsilon)$. Hence $C_{0,1}$ and $C_{0,1}'$ tend to be tangent when $\delta\to 0$, and we conclude again by Theorem \ref{existence min catenoids}.
		\ep
		
		In the rest of the section, we will denote by $\Cat_{n,j}$ the minimal catenoid whose asymptotic boundary is $C_{n,j}\cup C_{n,j}'$, and by $\SCat_{n,j}$ the closure in $\HH$ of the connected component of $\HH\setminus \Cat_{n,j}$ that contains the totally geodesic planes spanned by $C_{n,j}$ and $ C_{n,j}'$.

		\subsection{A version of the asymptotic Plateau problem}
		
		We now adapt the existence of solutions to the asymptotic Plateau problem developed in Section \ref{sec:APP} to construct minimal disks with asymptotic boundary the image of $f$ as above, contained in the complement of a suitably chosen family of solid catenoids $\Cat_{n,j}$.
		
		\begin{lem}\label{lemma Plateau catenoids}
			Let $\Lambda$ be the image of $f$, and let $\iota:\Z\to\{0,1\}$. There exists a stable minimal disk $\Sigma_\iota$ contained in $$\Omega_\iota:=\HH\setminus\bigcup_{n\in\Z}\SCat_{n,\iota(n)}$$ 
			with $\partial_\infty\Sigma_\iota=\Lambda$. 
		\end{lem}
		
\begin{figure}[htb]
			\centering
			\includegraphics[width=.8\textwidth]{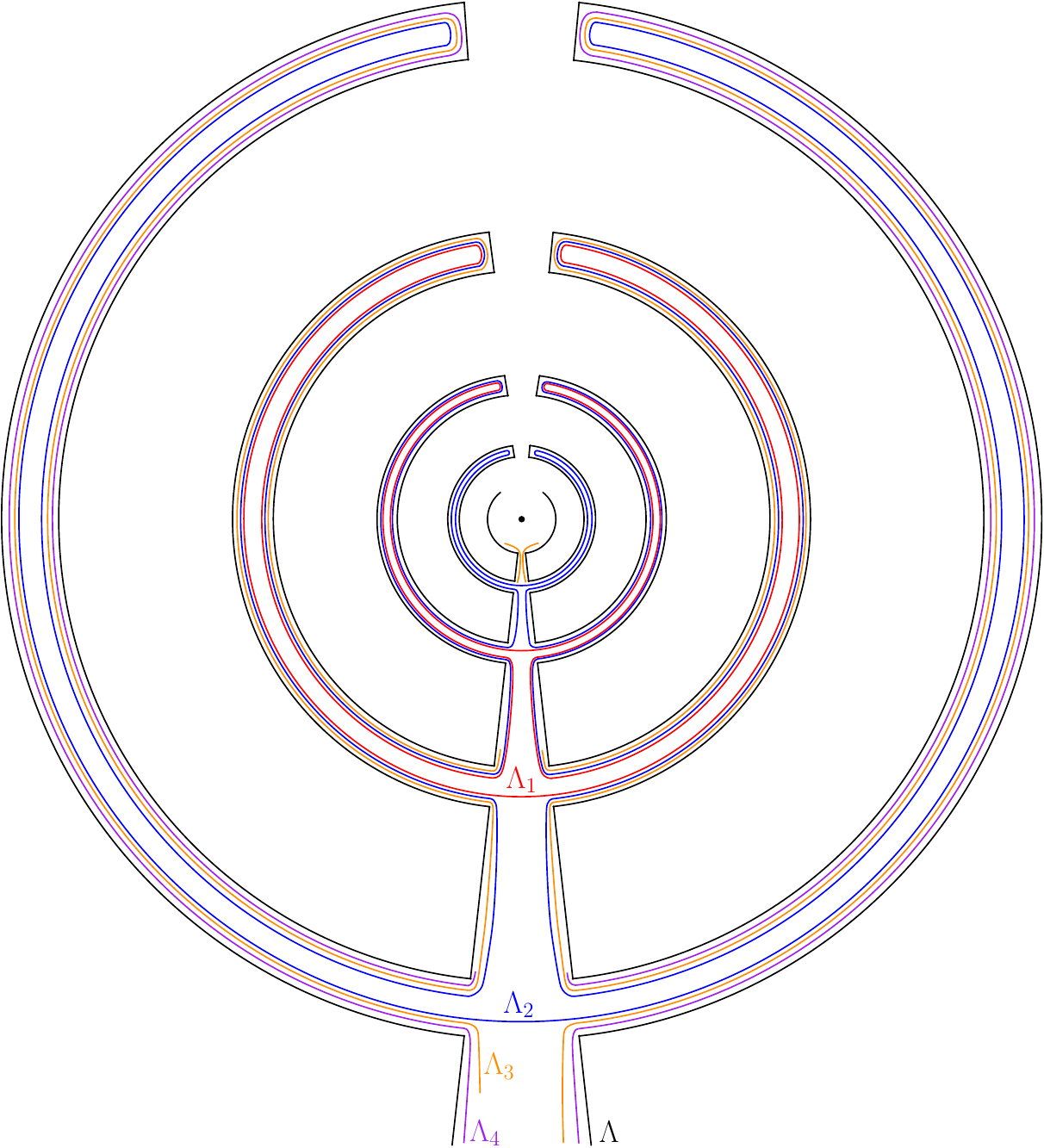} 
			\caption{\small The monotone sequence $\Lambda_n$ of smooth Jordan curves converging to $\Lambda$.}\label{fig:Jordancurve_approximate}
		\end{figure}

		\bp
		We follow the same strategy as in the proofs of Theorem \ref{andersonthm} and Lemma \ref{lemma Plateau}. First, let 
		$\Lambda_n$ be a monotone sequence of smooth Jordan curves disjoint from $\Lambda$, converging to $\Lambda$, and contained in the connected component of $\C\mathrm P^1\setminus\Lambda$ that contains the circles $C_{j,1}$ and $C'_{j,1}$, for all $j$. We require moreover that $\Lambda_n$ does not intersect any of the circles $C_{j,1}$ and $C'_{j,1}$, and that, among them, the connected component  of $\C\mathrm P^1\setminus\Lambda_n$ disjoint from $\Lambda$ contains the circles $C_{-n,1},C_{-n+1,1},\ldots,C_{n-1,1},C_{n,1}$ and $C'_{-n,1},C'_{-n+1,1},\ldots,C'_{n-1,1},C'_{n,1}$. It is very easy to see that such a sequence $\Lambda_n$ exists, see Figure \ref{fig:Jordancurve_approximate}.

		As we did for Lemma \ref{lemma Plateau}, we claim that, for each $\Lambda_n$, we can solve the asymptotic Plateau problem in $\Omega_\iota$, that is, we can find a properly embedded minimal surface $\Sigma_n$ contained in $\Omega_\iota$, area-minimizing in $\Omega_\iota$, with $\partial_\infty\Sigma_n=\Lambda_n$. Assuming this, the arguments of \cite[Lemma 3.2]{Cos08} show that  $\Sigma_n$ are a monotone sequence, and since each $\Sigma_n$ is area-minimizing in $\Omega$, it is stable, and thus we have $\|A_{\Sigma_n}\|\leq C$ for some $C>0$ and for all $n$. Hence we can apply Lemma \ref{lemma monotone convergence} and find a properly embedded stable minimal surface $\Sigma$, homeomorphic to a disk, which is the limit of the $\Sigma_n$. Using that every $\Sigma_n$ is contained in $\Omega_\iota$, the limit is contained in the closure of $\Omega_\iota$. However, by the strong maximum principle, if $\Sigma$ touched the boundary of $\Omega_\iota$, which has vanishing mean curvature, then it would coincide with a connected component of the boundary, which is a minimal catenoid $\Cat_{n,\iota(n)}$. But this is not possible because $\partial_\infty\Sigma=\Lambda$ is disjoint from all the circles $C_{n,\iota(n)}$ and $C'_{n,\iota(n)}$ by construction.
		
\begin{figure}[htbp]
			\centering
			\begin{minipage}[c]{.31\textwidth}
				\centering
				\includegraphics[height=4.5cm]{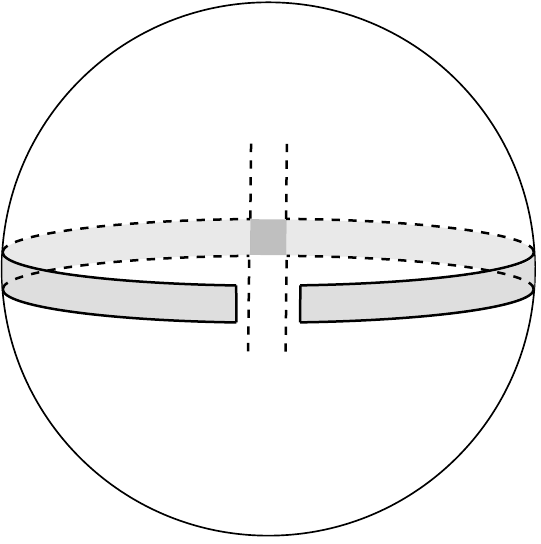} 
			\end{minipage}%
			\hspace{2mm}
			\begin{minipage}[c]{.31\textwidth}
				\centering
				\includegraphics[height=4.5cm]{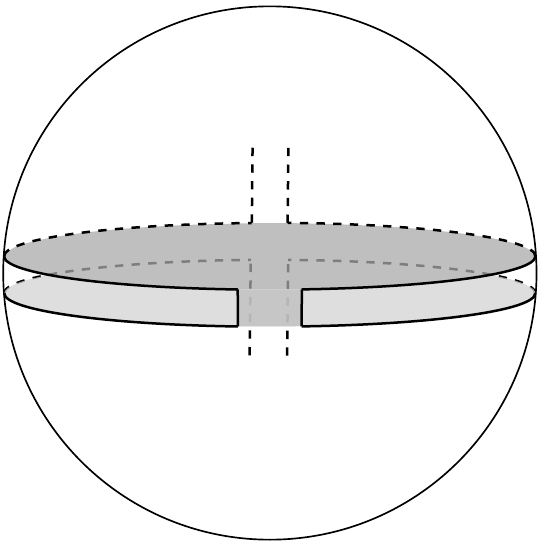}
			\end{minipage}
			\hspace{3mm}
			\begin{minipage}[c]{.31\textwidth}
				\centering
				\includegraphics[height=4.5cm]{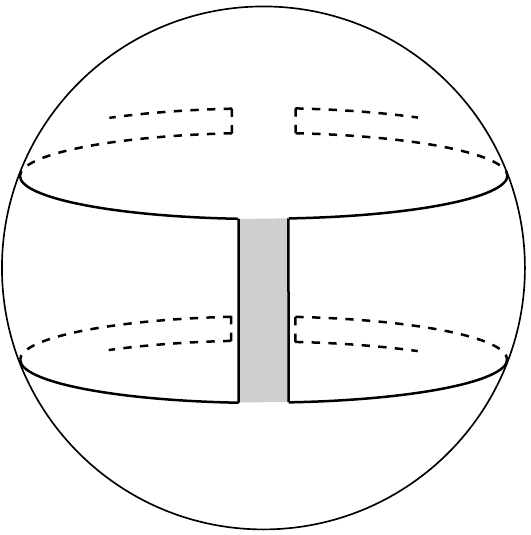}
			\end{minipage}
			\caption{The construction of the disk whose boundary is the curve $\Lambda_n$, which is represented here as piecewise smooth, although it is a smooth approximation of $\Lambda$. For any $j$ such that $|j|\leq n$, at the $j$-th level one chooses a disk that avoids the solid catenoid $\SCat_{j,0}$ if $\iota(j)=0$ (left figure), or a disk that avoids the solid catenoid $\SCat_{j,1}$ if $\iota(j)=1$ (middle figure). To connect the disks constructed for adjacent levels, one uses a band as in the right figure. Altogether, one obtains a disk whose boundary is $\Lambda_n$ that avoids the solid catenoids $\SCat_{-n,\iota(-n)}\ldots,\SCat_{n,\iota(n)}$. Clearly this disk can be arranged so as to also avoid the solid catenoids $\SCat_{j,0}$ and $\SCat_{j,1}$ for $|j|>n$. \label{fig:disks123new}}
		\end{figure}

		We thus need to show the claim. This is again an adaptation of the arguments of \cite{Som04}. To simplify the notation, let $\hat\Lambda$ denote one of the smooth Jordan curves $\Lambda_n$, contained in the boundary at infinity of $\Omega_\iota$. It is easy to see that $\hat\Lambda$ bounds a smooth disk in $\Omega_\iota$, that we call $\hat D$. See Figure \ref{fig:disks123new}. Now, we pick a sequence $\gamma_k$ of simple closed curves in $\hat D$ converging to $\hat\Lambda$. To be concrete, if $F:\mathbb D\to\hat D$ is a diffeomorphism, consider $\gamma_k:=F(\{|z|=1-1/n\})$. Using \cite[Section 1]{MY19} or \cite[Section 6]{HS88}, we pick the solutions of the finite Plateau problem for $\gamma_k$ inside the intersection of $\Omega_\iota$ with the half-spaces used by  Soma in \cite{Som04}. It is possible to apply \cite[Section 1]{MY19} or \cite[Section 6]{HS88} because each $\gamma_k$ is on $\hat D$, hence is nullhomotopic, and moreover (exactly as in Lemma \ref{lemma Plateau}) the intersection of $\Omega_\iota$ with a family of half-spaces is sufficiently convex, up to applying a small perturbation, since transversality of intersection is generic. The arguments of Soma then show that one can take the limit of these solutions of the finite Plateau problem for $\gamma_k$ to obtain a properly embedded minimal disk $\Sigma_n$ with asymptotic boundary $\hat\Lambda$. A priori $\Sigma_n$ is contained in the closure of $\Omega_\iota$, but by the strong maximum principle, it must actually be contained in $\Omega_\iota$.
		\ep

		\subsection{Conclusion of the proof}
		
		We are now ready to prove Theorem \ref{thm:uncountably}.
		
		\bp[Proof of Theorem \ref{thm:uncountably}]
		Let $\Lambda$ be the image of the continuous injective map $f$ as in Subsection \ref{subsec:construction jordan curve}, and fix $\iota:\Z\to\{0,1\}$. By Lemma \ref{lemma Plateau catenoids}, there exists a stable minimal disk $\Sigma_\iota$ contained in $$\Omega_\iota=\HH\setminus\bigcup_{n\in\Z}\SCat_{n,\iota(n)}$$ 
		with $\partial_\infty\Sigma_\iota=\Lambda$.

		It only remains to show that, if $\iota_0\neq \iota_1$, then $\Sigma_{\iota_0}\neq \Sigma_{\iota_1}$. By hypothesis, there exists $m\in\Z$ such that $\iota_0(m)\neq \iota_1(m)$. Up to switching the roles, suppose $\iota_0(m)=0$ and $\iota_1(m)=1$. Then, by construction, $\Sigma_{\iota_0}$ is contained in the complement of $\SCat_{m,0}$, whereas $\Sigma_{\iota_1}$ is contained in the complement of $\SCat_{m,1}$. But the solid catenoids $\SCat_{m,0}$ and $\SCat_{m,1}$ are linked in the sense of \cite[Section 3.4]{HW15}. Then \cite[Lemma 3.12]{HW15} implies that $\Lambda$ is not nullhomotopic in $(\HH\cup\partial_\infty\HH)\setminus(\SCat_{m,0}\cup\SCat_{m,1})$. Now, if $\Sigma_{\iota_0}= \Sigma_{\iota_1}=:\Sigma$, then $\Sigma$ would be a disk contained in the complement of $\SCat_{m,0}\cup\SCat_{m,1}$, hence contradicting that $\Lambda$ is not nullhomotopic. This concludes the proof.
		\ep

		\appendix 
		
		\section{Small curvatures}\label{appendix}
		
		In this appendix, we explain how to obtain Theorems \ref{exponential diffeo} and \ref{thm:weakly small is properly emb}. These results are well-known under the assumption that $\Sigma$ is properly embedded and has small curvatures (see for instance \cite{Eps84}). However, we need to relax the hypotheses here, for instance to assume that  $\Sigma$ has only weakly small curvatures. We follow here the approach of \cite[Section 4]{ES22}.
		
		The fundamental observation is that, if $\Sigma$ is a hypersurface in $\HHnp$ of weakly small curvatures and $\gamma:I\to\Sigma$ is a geodesic for the first fundamental form of $\Sigma$, then the acceleration of $\gamma$ is
		\begin{equation}\label{eq:acceleration}
			\nabla^h_{\gamma'(t)}\gamma'(t)=A(\gamma'(t),\gamma'(t))N(\gamma(t))~,
		\end{equation}
		where we recall that $\nabla^{h}$ is the Levi-Civita connection of $h$, the hyperbolic metric of $\HHnp$.
		Therefore, by the weakly small curvature assumption,
		\begin{equation}\label{eq:norm acceleration}
			\|\nabla^h_{\gamma'(t)}\gamma'(t)\|\leq 1~.
		\end{equation}
		
		This leads us to study curves satisfying the bound \eqref{eq:norm acceleration} on the acceleration. Let $S$ be a horosphere or a sphere in $\HHnp$. (By sphere, here we always mean the set of points at distance $r>0$ from a given point of $\HHnp$.) The complement of $S$ has two connected components, one of which is geodesically convex, and the other is not. We call the former the (open) convex side of $S$, and the latter the (open) concave side of $S$. The closed convex (resp. concave) side of $S$ is the union of $S$ with the open convex (resp. concave) side. The following lemma is an improvement of \cite[Lemma 4.9]{ES22}, following a similar strategy of proof.
		
		\bl\label{lemma:curve concave side}
		Let $\gamma:[a,b]\to\HHnp$ be a smooth curve satisfying \eqref{eq:norm acceleration}. Then 
		\begin{itemize}
			\item The image of $\gamma$ lies in the open concave side of any  sphere tangent to $\gamma$, except for the tangency point. 
			\item The image of $\gamma$ lies in the closed concave side of any horosphere tangent to $\gamma$. 
		\end{itemize}
		\el 
		\bp
		It suffices to prove the first item, as the second item follows by taking the limit of tangent spheres as the radius goes to infinity. 
		
		Assume that the tangency point is, in the upper half-space model $\gamma(0)=(0,\ldots,0,1)$, that  $\gamma$ is parameterized by arclength, that $\gamma'(0)=(1,0,\ldots,0)$ and that the tangent horosphere is $\{x_n=1\}$. Hence the tangent spheres $S_r$ of radius $r$ are, in a neighborhood of $\gamma(0)$, graphs of functions $f_r:U_r\to\R$, where $U_r\subset\R^n$ is a sufficiently small ball, and by symmetry $f_r$ only depends on the distance from the origin in $\R^n$. 
		
		First we show $\gamma$ lies on the open concave side of any tangent sphere $S_r$ for small $t$.
		We use a subscript $(\cdot)_{n+1}$ to denote the last coordinate in the upper half-space model. Thus $\gamma_{n+1}(0)=\gamma'_{n+1}(0)=0$.
		By a direct computation of the Christoffel symbol $\Gamma_{11}^n=1$, we see that
		$$(\nabla^h_{\gamma'}\gamma')_{n+1}(0)=\gamma_{n+1}''(0)+1~.$$
		Since $\gamma$ has small acceleration, we obtain 
		$$\gamma_{n+1}''(0)\leq |(\nabla^h_{\gamma'}\gamma')_{n+1}(0)|-1\leq \|(\nabla^h_{\gamma'}\gamma')(0)\|-1\leq 0~.$$
		Now, denote $\hat\gamma$ the projection of $\gamma$ to $\R^n$. By an easy computation following similar lines, $(f_r\circ\hat\gamma)(0)=(f_r\circ\hat\gamma)'(0)=0$ and $(f_r\circ\hat\gamma)''(0)=-1+1/\tanh(r)>0$. Hence for $t\in(0,\epsilon)$, $\gamma_{n+1}(t)<f_r(t)$ and this shows that $\gamma$ stays for small times in the open concave side of $S_r$.
		
		It remains to show that $\gamma(t)$ stays in in the open concave side of $S_r$ for every $t>0$. Suppose by contradiction that for some $t_0$, $\gamma(t_0)\in S_r$. Let $p$ be the center of $S_r$. Then the function $t\mapsto d(p,\gamma(t))$ is equal to $r$ for $t=0$ and $t=t_0$, and is larger than $r$ for $t\in (0,\epsilon)$. Hence it admits a maximum point $t_{\max}$. This implies that $\gamma$ is tangent to the  sphere centered at $p$ of radius $d(p,\gamma(t_{\max}))$ and thus contradicts the first part of the proof.
		\ep
		
		Having established Lemma \ref{lemma:curve concave side}, the following lemma follows immediately.
		
		\bl \label{lemma:surface concave side}
		Let $\Sigma\subset\HHnp$ be an immersed hypersurface of weakly small  curvatures, whose first fundamental form is geodesically convex. Then, for every $p\in\Sigma$, 
		\begin{itemize}
			\item $\Sigma\setminus\{p\}$ is contained in the open concave side of every  sphere tangent to $\Sigma$ at $p$.
			\item $\Sigma$ is contained in the closed concave side of every  horosphere tangent to $\Sigma$ at $p$.
		\end{itemize}
		\el
		\bp
		Let $p\in\Sigma$ and let $S$ be a tangent sphere. Given any other $q$, by geodesic convexity of $\Sigma$, we can find a geodesic segment $\gamma$ joining $p$ and $q$. As observed in \eqref{eq:norm acceleration}, $\gamma$ satisfies the hypothesis of Lemma \ref{lemma:curve concave side}, hence $q$ is contained in the open concave side of $S$. The proof of the second item is identical.
		\ep
		
		\begin{rem}
			Since every point of $\Sigma$ has a geodesically convex neighborhood, the conclusion of Lemma \ref{lemma:surface concave side} holds in a neighborhood of every point of $\Sigma$. It holds globally if $\Sigma$ is complete, by the Hopf-Rinow Theorem.
		\end{rem}
		
		Now, the proof of Theorem \ref{thm:weakly small is properly emb} follows exactly the proof of \cite[Proposition 4.15]{ES22}, observing that a properly immersed hypersurface has complete first fundamental form.

		To prove Theorem \ref{exponential diffeo}, we need the following well-known result, which follows from a direct computation (see for instance \cite[formula $(2.2)$]{HW13} or \cite[Lemma 4.5]{ES22}).
		
		\bl\label{lemma principal curvatures}
		Let $\Sigma$ be an embedded hypersurface in $\HHnp$ of weakly small curvatures, and let $N$ be its unit normal vector. Then for every $t\in\R$, the map $\iota_t:\Sigma\to\HHnp$ defined by $\iota_t(x)=\exp_x(tN(x))$ is an immersion, and its principal curvatures satisfy the identity
		\be\label{nfcurv}
		\lambda^t_i(\iota_t(x)) = {\frac{\lambda_i(x)-\tanh(t)}{1-\lambda_i(x)\tanh(t)}}~,
		\ene
		where $\lambda_1,\ldots,\lambda_n$ are the principal curvatures of $\Sigma$.
		\el
		\bp[Proof of Theorem \ref{exponential diffeo}]
		It follows from Lemma \ref{lemma principal curvatures} that the map $F(x,t)=\iota_t(x)$ has injective differential, hence is a local diffeomorphism. 
		
		To prove (i), we show that $F$ is injective if $\Sigma$ is geodesically convex. Suppose that $p=F(x_1,t_1)=F(x_2,t_2)$. Let $S_{t_i}$ be the geodesic spheres centered at $p$ of radius $t_i$, and observe that $S_{t_i}$ is tangent to $\Sigma$ at $x_i$. By Lemma \ref{lemma:surface concave side}, $\Sigma$ is contained, except for the tangency points, in the open concave side of both $S_{t_1}$ and $S_{t_2}$. This implies that $t_1=t_2$, because otherwise the open concave side of $S_{t_1}$ would intersect the convex side of $S_{t_2}$ (or vice versa), leading to a contradiction. By Lemma \ref{lemma:surface concave side} again, $\Sigma\setminus\{x_1\}$ is contained in the open concave side of $S_{t_1}=S_{t_2}$, and thus $x_2=x_1$. Hence $F$ is a diffeomorphism onto its image.
		
		When $\Sigma$ is in addition complete, as in (ii), we show that $F$ is surjective, and thus a diffeomorphism. Let $p$ be any point in $\HHnp$, let $r=d(p,\Sigma)$ and let $x\in\Sigma$ be a point that realizes this distance (that exists by completeness of $\Sigma$). Then the geodesic sphere of center $p$ and radius $r$ is tangent to $\Sigma$ at $x$, and $p=F(x,r)$. 
		
		It follows immediately that each $\iota_t$ is an embedding under the assumption of (i), and a proper embedding under the assumption of (ii). In both cases, it has weakly small curvatures by Lemma \ref{lemma principal curvatures}. 
		
		When $\Sigma$ is minimal, Lemma \ref{lemma principal curvatures} again implies the last statement on the sign of the mean curvature of $\iota_t$. Indeed, $\lambda_i^t(x)$ is a non-increasing function of $t$, and is decreasing unless $\lambda_i(x)=\pm 1$, in which case it is constant. So the mean curvature of $\iota_t$ at $x$ is a decreasing function of $t$, and becomes therefore negative for any $t>0$, unless all the principal curvatures are equal to $\pm 1$ --- equivalently, unless $\|A_\Sigma(x)\|^2=n$. But (see also \cite[Section 3]{HL21}) the subset $\{x\in\Sigma\,|\,\|A_\Sigma(x)\|^2=n\}$ has zero measure by Lemma \ref{zero measure} below, which concludes the proof. 
		\ep
		
		\begin{lem}\label{zero measure}
			Let $\Sigma$ a minimal immersion of weakly small curvature. The set $\{x\in\Sigma\,|\,\|A_\Sigma(x)\|^2=n\}$ has zero measure.
		\end{lem}
		\bp
		The subset where $\|A_\Sigma(x)\|^2=n$ is the zero set of the real analytic function $\|A_\Sigma\|^2-n$, hence (see \cite{Mit20}) it is either $\Sigma$ or a subset of zero measure. But it cannot be $\Sigma$ because in that case $\|A_\Sigma\|^2\equiv n$ implies that all principal curvatures are constantly equal to $\pm 1$, and this is not possible by Cartan's classification of isoparametric hypersurfaces in $\HHnp$ (\cite{Cartan}).
		\ep
		
		\begin{rem}
			When $n=2$, it is more immediate to see that the principal curvatures of a minimal surface cannot be identically $\pm 1$. Indeed in that case one has the formula $K_\Sigma=(1/4)\Delta(\log(\|A_\Sigma\|^2))$  for the curvature of the first fundamental form (see \cite[Lemma 3.11]{KS07}, \cite[Theorem 2.2]{ww20} or \cite[Appendix B]{EES22}). So if $\|A_\Sigma\|^2\equiv 2$ then $\Sigma$ should be intrinsically flat, but by Gauss' equation the curvature would be identically equal to $K_\Sigma=-1+\lambda_1\lambda_2=-2$.
		\end{rem}
		
		We have already observed that, for $n=2$, if the Jordan curve $\Lambda$ spans a surface of strongly small curvatures, then $\Lambda$ is a quasicircle, and this in turns implies that $\Lambda$ has finite width. We conclude this appendix by a direct proof, that works in any dimensions.
		
		\bl\label{lemma finite width}
		Let $\Lambda$ be a topologically embedded $(n-1)$-sphere $\Lambda$ in $\Sph_\infty^n$ and $\Sigma$ be a properly embedded hypersurface of strongly small curvatures such that $\partial_\infty\Sigma=\Lambda$. Then 
		$\Cal(\Lambda)$ has finite width.	\el
		\bp
		Assume that the principal curvatures $\lambda_1,\ldots,\lambda_n$ of $\Sigma$ are smaller than $1-\epsilon$ in absolute value. By Lemma \ref{nfcurv}, for $t_0 > \tanh^{-1}(1-\epsilon)$ the hypersurface $\Sigma_{t_0}$ (resp. $\Sigma_{-t_0}$) has negative (resp. positive) principal curvatures, that is, it is locally concave (resp. convex) with respect to the direction of the normal vector $N$. Since $\Sigma_{\pm t_0}$ are properly embedded, they bound a geodesically convex region $U$, which thus contains the $\Cal(\Lambda)$. By Theorem \ref{exponential diffeo}, every point in $\Cal(\Lambda)$ lies on a geodesic segment of the form $t\mapsto F(x,t)$, for $t\in [-t_0,t_0]$. Hence for every $x\in \Cal(\Lambda)$, $d(x,\partial^{+}\Cal(\Lambda))+d(x,\partial^{-}\Cal(\Lambda))\leq 2t_0$. This shows that $w(\Lambda)\leq 2t_0$.  
		\ep

	\bibliographystyle{amsbook}

\begin{thebibliography}{99}
		\bibitem[Ahl66]{Ahl66}
		Lars~V. Ahlfors, \emph{Lectures on quasiconformal mappings}, Princeton,
		{N}.{J}.-{Toronto}-{New} {York}-{London}: {D}. {Van} {Nostrand} {Company}.
		{Inc}. 146 p., 1966.
		
		\bibitem[AM10]{AM10}
		Spyridon Alexakis and Rafe Mazzeo, \emph{Renormalized area and properly
			embedded minimal surfaces in hyperbolic 3-manifolds}, Comm. Math. Phys.
		\textbf{297} (2010), no.~3, 621--651.
		
		\bibitem[And83]{And83}
		Michael Anderson, \emph{Complete minimal hypersurfaces in hyperbolic
			{$n$}-manifolds}, Comment. Math. Helv. \textbf{58} (1983), no.~2, 264--290.
		
		\bibitem[Bac24]{Ba24}
		Aidan Backus, \emph{Minimal laminations and level sets of 1-harmonic
			functions}, ArXiv 2311.01541, to appear in Journal of Geometric Analysis
		(2024).
		
		\bibitem[BDMS21]{BDMS21}
		Francesco Bonsante, Jeffrey Danciger, Sara Maloni, and Jean-Marc Schlenker,
		\emph{Quasicircles and width of {J}ordan curves in {$\Bbb C\Bbb P^1$}}, Bull.
		Lond. Math. Soc. \textbf{53} (2021), no.~2, 507--523.
		
		\bibitem[Bis19]{Bis19}
		Christopher~J. Bishop, \emph{Weil-{P}etersson curves, conformal energies,
			$\beta$-numbers, and minimal surfaces},
		https://www.math.stonybrook.edu/~bishop/papers/wpbeta.pdf (2019).
		
		\bibitem[Bro61]{monotoneunion}
		Morton Brown, \emph{The monotone union of open {{\(n\)}}-cells is an open
			{{\(n\)}}-cell}, Proc. Am. Math. Soc. \textbf{12} (1961), 812--814.
		
		\bibitem[Bro23]{b23}
		Samuel Bronstein, \emph{Almost-{F}uchsian structures on disk bundles over a
			surface}, ArXiv:2305.06665 (2023).
		
		\bibitem[BT16]{BT16}
		Jacob Bernstein and Giuseppe Tinaglia, \emph{Topological type of limit
			laminations of embedded minimal disks}, J. Differ. Geom. \textbf{102} (2016),
		no.~1, 1--23.
		
		\bibitem[Car38]{Cartan}
		{\'E}lie Cartan, \emph{Familles de surfaces isoparametriques dans les espaces
			{\`a} courbure constante}, Ann. Mat. Pura Appl. (4) \textbf{17} (1938),
		177--191 (French).
		
		\bibitem[CF19]{CF19}
		Daryl Cooper and David Futer, \emph{Ubiquitous quasi-{F}uchsian surfaces in
			cusped hyperbolic 3-manifolds}, Geom. Topol. \textbf{23} (2019), no.~1,
		241--298.
		
		\bibitem[CM04]{CMIV}
		Tobias~H. Colding and William P.~II Minicozzi, \emph{The space of embedded
			minimal surfaces of fixed genus in a 3-manifold. {IV}: {Locally} simply
			connected}, Ann. Math. (2) \textbf{160} (2004), no.~2, 573--615.
		
		\bibitem[CMN22]{CMN22}
		Danny Calegari, Fernando~C. Marques, and Andr{\'e} Neves, \emph{Counting
			minimal surfaces in negatively curved 3-manifolds}, Duke Math. J.
		\textbf{171} (2022), no.~8, 1615--1648.
		
		\bibitem[Cos08]{Cos08}
		Baris Coskunuzer, \emph{Properly embedded least area planes in {Gromov}
			hyperbolic 3-spaces}, Proc. Am. Math. Soc. \textbf{136} (2008), no.~4,
		1427--1432.
		
		\bibitem[Cos11]{c11genericuniqueness}
		\bysame, \emph{On the number of solutions to the asymptotic {P}lateau problem},
		J. G\"{o}kova Geom. Topol. GGT \textbf{5} (2011), 1--19.
		
		\bibitem[Cos14]{co14}
		\bysame, \emph{Asymptotic {P}lateau problem: a survey}, Proceedings of the
		{G}\"{o}kova {G}eometry-{T}opology {C}onference 2013, G\"{o}kova, 2014,
		pp.~120--146.
		
		\bibitem[Cos16]{Cos16}
		\bysame, \emph{Asymptotic {$H$}-{P}lateau problem in {$\Bbb H^3$}}, Geom.
		Topol. \textbf{20} (2016), no.~1, 613--627.
		
		\bibitem[Cos19]{coskunuzer2019nestedsequences}
		\bysame, \emph{Embedded H-planes in hyperbolic 3-space}, Transactions of
		the American Mathematical Society \textbf{371} (2019), no.~2, 1253--1269.
		
		\bibitem[dMG87]{MG87}
		Jonas de~Miranda~Gomes, \emph{Spherical surfaces with constant mean curvature
			in hyperbolic space}, Bol. Soc. Bras. Mat. \textbf{18} (1987), no.~2, 49--73.
		
		\bibitem[EES22a]{ES22}
		Christian El~Emam and Andrea Seppi, \emph{On the {G}auss map of equivariant
			immersions in hyperbolic space}, J. Topol. \textbf{15} (2022), no.~1,
		238--301.
		
		\bibitem[EES22b]{EES22}
		\bysame, \emph{Rigidity of minimal {Lagrangian} diffeomorphisms between
			spherical cone surfaces}, J. {\'E}c. Polytech., Math. \textbf{9} (2022),
		581--600.
		
		\bibitem[Eps84]{Eps84}
		Charles~L. Epstein, \emph{Envelopes of horospheres and {W}eingarten surfaces in
			hyperbolic 3-space}, unpublished manuscript.
		
		\bibitem[Eps86]{Eps86}
		\bysame, \emph{The hyperbolic {G}auss map and quasiconformal reflections}, J.
		Reine Angew. Math. \textbf{372} (1986), 96--135.
		
		\bibitem[Fin21]{f21}
		Joel Fine, \emph{Knots, minimal surfaces and j-holomorphic curves},
		arxiv:2112.07713 (2021).
		
		\bibitem[Gab97]{Gab97}
		David Gabai, \emph{On the geometric and topological rigidity of hyperbolic
			{$3$}-manifolds}, J. Amer. Math. Soc. \textbf{10} (1997), no.~1, 37--74.
		
		\bibitem[GHW10]{GHW10}
		Ren Guo, Zheng Huang, and Biao Wang, \emph{Quasi-{F}uchsian three-manifolds and
			metrics on {\TS}}, Asian J. Math. \textbf{14} (2010), no.~2, 243--256.
		
		\bibitem[GLT22]{glt88}
		Michael Gromov, Herbert~Blaine Lawson, Jr., and William Thurston,
		\emph{Hyperbolic 4-manifolds and conformally flat 3-manifolds}, Collected
		works of {W}illiam {P}. {T}hurston with commentary. {V}ol. {I}. {F}oliations,
		surfaces and differential geometry, Amer. Math. Soc., Providence, RI, 2022,
		Reprint of [1001446], pp.~667--685.
		
		\bibitem[GLVP21]{GLP21}
		Marco Guaraco, Vanderson Lima, and Franco Vargas~Pallete, \emph{Mean curvature
			flow in homology and foliations of hyperbolic $3$-manifolds}, 2021.
		
		\bibitem[Gov]{mathstack}
		Dejan Govc, \emph{Connectedness of the boundary}, Mathematics Stack Exchange,
		URL:https://math.stackexchange.com/q/170355 (version: 2013-12-27).
		
		\bibitem[GS92]{GS92}
		Frederick~P. Gardiner and Dennis~P. Sullivan, \emph{Symmetric structures on a
			closed curve}, Am. J. Math. \textbf{114} (1992), no.~4, 683--736.
		
		\bibitem[GS00]{GS00}
		Bo~Guan and Joel Spruck, \emph{Hypersurfaces of constant mean curvature in
			hyperbolic space with prescribed asymptotic boundary at infinity}, Amer. J.
		Math. \textbf{122} (2000), no.~5, 1039--1060.
		
		\bibitem[Has15]{has15}
		Joel Hass, \emph{Minimal fibrations and foliations of hyperbolic
			$3$-manifolds}, Preprint (2015).
		
		\bibitem[HL21]{HL21}
		Zheng Huang and Ben Lowe, \emph{Beyond almost {F}uchsian space}, ArXiv
		2104.11284 (2021).
		
		\bibitem[HS88]{HS88}
		Joel Hass and Peter Scott, \emph{The existence of least area surfaces in
			3-manifolds}, Trans. Am. Math. Soc. \textbf{310} (1988), no.~1, 87--114.
		
		\bibitem[HW13]{HW13}
		Zheng Huang and Biao Wang, \emph{On almost-{F}uchsian manifolds}, Trans. Amer.
		Math. Soc. \textbf{365} (2013), no.~9, 4679--4698.
		
		\bibitem[HW15]{HW15}
		\bysame, \emph{Counting minimal surfaces in quasi-{F}uchsian manifolds}, Trans.
		Amer. Math. Soc. \textbf{367} (2015), 6063--6083.
		
		\bibitem[HW19]{HW19}
		\bysame, \emph{Complex length of short curves and minimal foliations in closed
			hyperbolic three-manifolds fibering over the circle}, Proc. London Math. Soc.
		\textbf{118} (2019), no.~3, 1305--1327.
		
		\bibitem[Jia21]{j21}
		Ruojing Jiang, \emph{Counting essential minimal surfaces in closed negatively
			curved n-manifolds}, arxiv:2108.01796 (2021).
		
		\bibitem[JT03]{jt03}
		Luqu\'{e}sio~P. Jorge and Friedrich Tomi, \emph{The barrier principle for
			minimal submanifolds of arbitrary codimension}, Ann. Global Anal. Geom.
		\textbf{24} (2003), no.~3, 261--267.
		
		\bibitem[KM12]{KM12}
		Jeremy Kahn and Vladimir Markovic, \emph{Immersing almost geodesic surfaces in
			a closed hyperbolic $3$-manifold}, Ann. of Math. (2) \textbf{175} (2012),
		no.~3, 1127--1190.
		
		\bibitem[KMS23]{KMS23}
		Jeremy Kahn, Vladimir Markovic, and Ilia Smilga, \emph{Geometrically and
			topologically random surfaces in a closed hyperbolic three manifold}, ArXiv
		2309.02847 (2023).
		
		\bibitem[KS07]{KS07}
		Kirill Krasnov and Jean-Marc Schlenker, \emph{Minimal surfaces and particles in
			3-manifolds}, Geom. Dedicata \textbf{126} (2007), 187--254.
		
		\bibitem[KW21]{kw21}
		Jeremy Kahn and Alex Wright, \emph{Nearly {F}uchsian surface subgroups of
			finite covolume {K}leinian groups}, Duke Math. J. \textbf{170} (2021), no.~3,
		503--573.
		
		\bibitem[Lei06]{Lei06}
		Christopher~J. Leininger, \emph{Small curvature surfaces in hyperbolic
			3-manifolds}, J. Knot Theory Ramifications \textbf{15} (2006), no.~3,
		379--411.
		
		\bibitem[Low21]{Low21}
		Ben Lowe, \emph{Deformations of totally geodesic foliations and minimal
			surfaces in negatively curved 3-manifolds}, Geom. Funct. Anal. \textbf{31}
		(2021), no.~4, 895--929.
		
		\bibitem[Mit20]{Mit20}
		Boris~S. Mityagin, \emph{The zero set of a real analytic function}, Math. Notes
		\textbf{107} (2020), no.~3, 529--530.
		
		\bibitem[MR06]{MRclosure}
		William H.~III Meeks and Harold Rosenberg, \emph{The minimal lamination closure
			theorem}, Duke Math. J. \textbf{133} (2006), no.~3, 467--497.
		
		\bibitem[MY19]{MY19}
		William~W. Meeks, III and Shing-Tung Yau, \emph{The existence of embedded
			minimal surfaces and the problem of uniqueness}, Selected works of
		{S}hing-{T}ung {Y}au. {P}art 1, {V}ol. 1., Int. Press, Boston, MA, 2019,
		pp.~341--358.
		
		\bibitem[Ros06]{Ros06}
		Antonio Ros, \emph{One-sided complete stable minimal surfaces}, J. Differential
		Geom. \textbf{74} (2006), no.~1, 69--92.
		
		\bibitem[RST10]{RST10}
		Harold Rosenberg, Rabah Souam, and Eric Toubiana, \emph{General curvature
			estimates for stable {$H$}-surfaces in 3-manifolds and applications}, J.
		Differential Geom. \textbf{84} (2010), no.~3, 623--648.
		
		\bibitem[Rub05]{Rub05}
		J.~Hyam Rubinstein, \emph{Minimal surfaces in geometric 3-manifolds}, Global
		theory of minimal surfaces, Clay Math. Proc., vol.~2, AMS, 2005,
		pp.~725--746.
		
		\bibitem[Sam69]{Sam69}
		Hans Samelson, \emph{Orientability of hypersurfaces in $\mathbb{R}^n$}, Proc.
		Am. Math. Soc. \textbf{22} (1969), 301--302.
		
		\bibitem[San18]{San18}
		Andrew Sanders, \emph{Entropy, minimal surfaces and negatively curved
			manifolds}, Ergodic Theory Dynam. Systems \textbf{38} (2018), no.~1,
		336--370.
		
		\bibitem[Sep16]{Sep16}
		Andrea Seppi, \emph{Minimal discs in hyperbolic space bounded by a quasicircle
			at infinity}, Comment. Math. Helv. \textbf{91} (2016), no.~4, 807--839.
		
		\bibitem[Smi07]{Smi07}
		Graham Smith, \emph{An {Arzela}-{Ascoli} theorem for immersed submanifolds},
		Ann. Fac. Sci. Toulouse, Math. (6) \textbf{16} (2007), no.~4, 817--866.
		
		\bibitem[Som04]{Som04}
		Teruhiko Soma, \emph{Existence of least area planes in hyperbolic 3-space with
			co-compact metric}, Topology \textbf{43} (2004), no.~3, 705--716.
		
		\bibitem[Thu86]{Thu86}
		William~P. Thurston, \emph{Hyperbolic structures on 3-manifolds, {II}: Surface
			groups and 3-manifolds which fiber over the circle}, ArXiv 1998 (1986).
		
		\bibitem[Uhl83]{Uhl83}
		Karen~K. Uhlenbeck, \emph{Closed minimal surfaces in hyperbolic
			{$3$}-manifolds}, Seminar on minimal submanifolds, Ann. of Math. Stud., vol.
		103, Princeton Univ. Press, Princeton, NJ, 1983, pp.~147--168.
		
		\bibitem[Whi10]{Whi10}
		Brian White, \emph{The maximum principle for minimal varieties of arbitrary
			codimension}, Comm. Anal. Geom. \textbf{18} (2010), no.~3, 421--432.
		
		\bibitem[WW20]{ww20}
		Michael Wolf and Yunhui Wu, \emph{Non-existence of geometric minimal foliations
			in hyperbolic three-manifolds}, Comment. Math. Helv. \textbf{95} (2020),
		no.~1, 167--182.
		
	\end{thebibliography}

	\end{document}